\documentclass[preprint,12pt]{elsarticle}
\renewcommand{\u}{\mathbf{u}}
\newcommand{\halb}{\frac{1}{2}}

\newcommand{\w}{\mathbf{w}}
\newcommand{\q}{\mathbf{q}}
\newcommand{\be}{\begin{equation}}
\newcommand{\ee}{\end{equation}}
\newcommand{\bdm}{\begin{displaymath}}
\newcommand{\edm}{\end{displaymath}}
\newcommand{\bea}{\begin{eqnarray}}
\newcommand{\eea}{\end{eqnarray}}

\newcommand{\ar}{\phi_1\rho_1}
\newcommand{\arr}{\phi_2\rho_2}

\newcommand{\ub}{\textbf{v}_1}

\newcommand{\tn}{t^n}
\newcommand{\tnext}{t^{n+1}}
\newcommand{\xL}{x_{i-\frac{1}{2}}}
\newcommand{\xR}{x_{i+\frac{1}{2}}}
\newcommand{\yL}{y_{j-\frac{1}{2}}}
\newcommand{\yR}{y_{j+\frac{1}{2}}}
\newcommand{\zL}{z_{k-\frac{1}{2}}}
\newcommand{\zR}{z_{k+\frac{1}{2}}}

\usepackage{amsthm}
\usepackage{amsfonts}
\usepackage{latexsym}
\usepackage{wasysym}

\begin{document}

\begin{frontmatter}



\title{High Order Space-Time Adaptive WENO Finite Volume Schemes for Non-Conservative Hyperbolic Systems} 


\author[1]{Michael Dumbser\corref{corr1}\fnref{label2}}
\author[2]{Arturo Hidalgo}
\author[1]{Olindo Zanotti}
\address[1]{Laboratory of Applied Mathematics, Department of Civil, Environmental and Mechanical Engineering, 
          					   University of Trento,
          						 I-38123 Trento, Italy}
\address[2]{Departamento de Matem\'atica Aplicada y M\'etodos Inform\'aticos,
          							Universidad Polit\'ecnica de Madrid, 
          							28003-Madrid, Spain}
												

\fntext[label2]{\tt michael.dumbser@unitn.it (M.~Dumbser)}         


\begin{abstract}
We present a class of high order finite volume schemes for the solution of non-conservative hyperbolic systems that combines the one-step ADER-WENO finite volume approach with 
space-time adaptive mesh refinement (AMR). The resulting algorithm, which is particularly well suited for the treatment of material interfaces in compressible multi-phase flows, 
is based on: (i) high order of accuracy in space obtained through WENO reconstruction, (ii) a high order one-step time discretization via a local space-time discontinuous Galerkin 
predictor method, and (iii) the use of a path conservative scheme for handling the non-conservative terms of the equations. 
The AMR property with \textit{time accurate local time stepping}, which has been treated according to a \textit{cell-by-cell} strategy, strongly relies on the high order one-step time 
discretization, which naturally allows a high order accurate and consistent computation of the jump terms at interfaces between elements using different time steps. The new scheme 
has been succesfully validated on some test problems for the Baer-Nunziato model of compressible multiphase flows. 
\end{abstract}


\end{frontmatter}


\section{Introduction}
\label{Intro}

Physical phenomena governed by non-conservative hyperbolic systems arise in many scientific and technological areas, such as aerospace and automotive industry, geophysical flows, compressible multi-phase flows, 
oil and gas extraction, among many others. The mathematical modelling of these phenomena is very complicated and so far no universally accepted model exists. A common problem in many mathematical models used in the 
fields listed above is that the governing PDE system can not be written in fully conservative form. A special problem that concerns compressible multi-phase flows is the accurate resolution of material interfaces 
over long times. Among possible strategies to achieve a good resolution of material interfaces are the use of Lagrangian methods \cite{Maire2011,MaireMM1,MaireMM2,Despres2009,chengshu1,chengshu2,DumbserBoscheri}, 
ghost-fluid and level-set methods \cite{FedkiwEtAl1,FedkiwEtAl2,levelset1,levelset2,FerrariLevelSet}, little dissipative Riemann solvers combined with high order schemes \cite{TitarevTVD1,TitarevTVD2,OsherNC} and, 
last but not least, the use of adaptive mesh refinement (AMR). Therefore, in this paper we suggest to combine high order accurate WENO finite volume schemes with AMR for the solution of compressible multi-phase 
flow problems to assure an accurate resolution of the material interfaces. 

One particular and rather widespread mathematical model used in multi-phase flow applications is based on the  deflagration-to-detonation transition (DDT) theory, from which the Baer-Nunziato equations 
\cite{BaerNunziato1986} can be derived. \textit{Reduced} models have been proposed in \cite{Kapila2001,MurroneGuillard,ZeinHantke}, where the original Baer-Nunziato system has been simplified by carrying out an 
asymptotic analysis in the limit of stiff relaxation source terms.  

Other models in the field of multi-phase flows are those proposed and analyzed by Saurel and Abgrall \cite{SaurelAbgrall,AbgrallSaurel,SaurelCavitation,SaurelVapor} for the mixture of two compressible fluids,  
the depth-averaged debris-flow model by Pitman and Le \cite{PitmanLe} as well as single and multi-layer shallow water equations \cite{Castro2006,Pares2006,AbgrallKarni3}. 

All the models above can be cast into the following general form of a nonlinear system of PDE in multiple space dimensions 
\begin{equation}	\label{NCsyst}
	\frac{\partial \u}{\partial t}+\nabla\cdot\bf F(\u)+\bf B(\u)\cdot\nabla \u=\bf S(\u)\,,
\end{equation}
where $\u$ is the state vector;
${\bf F}(\u)=[\textbf{f}(\u),\textbf{g}(\u),\textbf{h}(\u)]$ is the flux tensor for the conservative part of the PDE system,
with $\textbf{f}(\u)$, $\textbf{g}(\u)$ and $\textbf{h}(\u)$ expressing the fluxes along the 
$x$, $y$ and $z$ directions, respectively;
${\bf B}(\u)=[\textbf{B}_1(\u),\textbf{B}_2(\u),\textbf{B}_3(\u)]$ represents the non-conservative part of the system, written in 
block-matrix notation. Finally, $\bf S(\u)$ is the source term, which may in principle be stiff.
When written in quasilinear form, the system (\ref{NCsyst}) becomes
\begin{equation}	\label{Csyst}
\frac{\partial \u}{\partial t}+ \bf{A}(\u) \cdot\nabla \u = \textbf{S}(\u)\,,
\end{equation}
where the matrix ${\bf A}(\u)=[{\bf A}_1,{\bf A}_2,{\bf A}_3]=\partial {\bf F}(\u)/\partial \u+{\bf B}(\u)$ contains both the conservative and the non-conservative contributions. 

Recent work on numerical schemes for systems of equations involving non-conservative terms, like Eq.~(\ref{NCsyst}), 
includes the family of so--called path-conservative schemes \cite{Castro2006,Castro2007,Castro2008,Pares2004} which are based on the theory proposed by Dal Maso, Le Floch and Murat \cite{DLMtheory} 
and are a generalization of the usual concept of conservative methods for systems of conservation laws. Note that the weak formulation of the Roe method by Toumi \cite{Toumi1992} can
also be considered as a path-conservative scheme. It has to be clearly stressed that path-conservative schemes have known deficiencies, which have been studied in detail in \cite{AbgrallKarni,NCproblems}. 

In this paper we provide the first implementation of high order path-conservative schemes for non-conservative systems of the type (\ref{NCsyst}) using the ADER approach 
together with space-time \textit{Adaptive Mesh Refinement} (AMR). In this respect, the present work can be considered an extension of the method proposed in \cite{Dumbser2012a}, where an 
ADER-WENO AMR scheme was proposed for the conservative case. 
ADER schemes, originally developed by E.F. Toro and collaborators in \cite{toro1} and extensively used in the context of hyperbolic problems 
\cite{titarevtoro,toro3,DumbserKaeser07,Balsara2009,Balsara2013}, are a class of numerical methods that obtain high order one-step schemes in time without the 
use of backward time levels, like in Adams-Bashforth type time integrators, and also without the use of substages, as used inside Runge-Kutta time integrators. The feature of high order one-step time  
integration will be actually the key for the construction of reasonably simple high order accurate AMR schemes together with time accurate local time stepping (LTS). For previous works including LTS
see, e.g., \cite{Flaherty,dumbserkaeser06d,TaubeMaxwell,stedg1,stedg2,DomainDecomp,CastroLTS,KrivodonovaLTS,Dumbser2012a}. 

However, the original ADER approach as proposed in \cite{toro1,titarevtoro,toro3} suffers from the drawback that it uses Taylor expansions in time where time derivatives are then replaced by 
spatial derivatives through a repeated use of the governing system of equations. This procedure, also known as the Cauchy-Kowalewski procedure, becomes rather cumbersome when dealing with complex systems 
of equations, and it fails completely in the presence of stiff source terms. A successful alternative to this strategy, that we also follow in this work, was proposed in \cite{DumbserEnauxToro}. There, 
a different formulation of the ADER approach was developed, where the Taylor series expansions and the Cauchy-Kovalewski procedure have been replaced by a local space-time Galerkin method, i.e. by a weak 
formulation of the PDE in space-time, see also the references \cite{Dumbser2008,DumbserZanotti,Balsara2009,Balsara2013} where this new version of the ADER approach has been used. 

As for the AMR aspects of our work, we have followed a {\em cell-by-cell} refinement strategy, which is particularly convenient within our high order one-step finite volume approach. 
Various examples of AMR schemes can be found in literature, for an overview see the original work by Berger et al. \cite{Berger-Oliger1984,Berger-Jameson1985,Berger-Colella1989,Berger-Leveque1998,Bell1994} 
and other implementations in  \cite{Mulet1,Balsara2001,KeppensAMR,FalleAMR1,FalleAMR2,Quirk1996,Carroll-Nellenback2011,Teyssier2002,Kurganov2005,Ziegler2008,Crash}. 

The outline of the paper is as follows: In Sect. \ref{sec:Numscheme} we present our numerical method, with some emphasis on the implementation of path-conservative schemes within the ADER approach. 
Sect. \ref{sec:AMR} is devoted to the description of the Adaptive Mesh Refinement infrastructure, while in Sect. \ref{sec:BN} the Baer-Nunziato equations are recalled. Sect. \ref{sec:tests} contains a 
set of numerical test problems and computational results to validate the proposed high order path-conservative ADER-WENO AMR schemes, and, finally, Sect. \ref{sec:concl} reports some conclusions of 
our work and possible future extensions.

%
\section{Numerical method}
\label{sec:Numscheme}

%
%
%
%
%
%

In order to obtain a numerical solution of the problem (\ref{NCsyst}) we use higher order finite volume methods in the context of the ADER framework. 
To simplify the presentation, we first describe all details of the algorithm for a uniform Cartesian grid. The AMR technique will be described later.  
We recall that within the finite volume methodology, the numerical solution of the evolved quantities is represented at the beginning of 
each time-step by piecewise constant cell averages. The update of these data and the computation of the corresponding numerical fluxes can be  
performed using higher order piecewise polynomials of degree $M$ that have to be reconstructed, starting from the underlying piecewise constant 
cell averages, see \cite{eno,shu_efficient_weno,TitarevToroWENO3D}. 

\subsection{The Finite Volume scheme for nonconservative systems}
\label{sec:FVM}

The system (\ref{NCsyst}) is written in Cartesian coordinates and in three space dimensions as
\begin{equation}	\label{NCsyst_cartesian}
	\frac{\partial \u}{\partial t}+\frac{\partial \textbf{f}}{\partial x}+\frac{\partial \textbf{g}}{\partial y}+\frac{\partial \textbf{h}}{\partial z}+\textbf{B}_1\frac{\partial \u}{\partial x}+\textbf{B}_2\frac{\partial \u}{\partial y}+\textbf{B}_3\frac{\partial\u}{\partial z}=\textbf{S}(\u)\,.
\end{equation}
In order to obtain the finite volume representation of (\ref{NCsyst_cartesian}), we discretize the computational domain $\Omega$ in space-time control volumes defined as ${\mathcal I}_{ijk}=I_{ijk}\times [t^n,t^n+\Delta t]=[\xL,\xR]\times[\yL,\yR]\times[\zL,\zR]\times [\tn,\tn+\Delta t]$, with $\Delta x_i=\xR-\xL$, $\Delta y_j=\yR-\yL$, $\Delta z_k=\zR-\zL$ and $\Delta t=\tnext-\tn$. Each space control volume $I_{ijk}$ defines a computational cell, which will be denoted so forth as ${\mathcal C}_m$, identified by its mono-index $m$, with $1<m<N_e$,  where $N_e$ is the number of computational cells in the domain. 
After integration of (\ref{NCsyst_cartesian}) over a space-time control volume ${\mathcal I}_{ijk}$ one obtains the following finite volume formulation: 
\begin{eqnarray} \label{FVformula}
\label{eq:finite_vol}
{\bar \u}_{ijk}^{n+1}&=&{\bar \u}_{ijk}^{n}-
\frac{\Delta t}{\Delta x_i}\left[\left(\textbf{f}_{i+\halb,j,k} \,\, -\textbf{f}_{i-\halb,j,k}\right)+\frac{1}{2} \,\,
\left({{D}}^x_{i+\halb,j,k} +{{D}}^x_{i-\halb,j,k}\right)\right]\nonumber\\
&& \hspace{9mm} -\frac{\Delta t}{\Delta y_j}\left[\left(\textbf{g}_{i,j+\halb,k}-\textbf{g}_{i,j-\halb,k}\right)+\frac{1}{2}
\left({{D}}^y_{i,j+\halb,k}+{{D}}^y_{i,j-\halb,k}\right)\right]\nonumber\\
&& \hspace{9mm} -\frac{\Delta t}{\Delta z_k}\left[\left(\textbf{h}_{i,j,k+\halb}-\textbf{h}_{i,j,k-\halb}\right)+\frac{1}{2}
\left({{D}}^z_{i,j,k+\halb}+{{D}}^z_{i,j,k-\halb}\right)\right]\nonumber\\
&& \hspace{9mm} + \Delta t({\bf \bar{S}}_{ijk}- {\bf \bar{P}}_{ijk})\,,
\end{eqnarray}

%
%
where
\begin{equation}
{\bar \u}_{ijk}^{n}=\frac{1}{\Delta x_i}\frac{1}{\Delta y_j}\frac{1}{\Delta z_k}\int_{x_{i-\halb}}^{x_{i+\halb}}\int_{y_{j-\halb}}^{y_{j+\halb}}\int_{z_{k-\halb}}^{z_{k+\halb}}{\u}(x,y,z,t^n)dz\,dy\,\,dx
\end{equation}
is the spatial average of the solution in the element ${ I}_{ijk}$ at time $\tn$, while
\begin{equation}
\label{averF}
{\bf f}_{i+\halb,jk}= \frac{1}{\Delta t}\frac{1}{\Delta y_j}\frac{1}{\Delta z_k} \hspace{-1mm}  \int \limits_{t^n}^{t^{n+1}} \! \int \limits_{y_{j-\halb}}^{y_{j+\halb}} \! \int \limits_{z_{k-\halb}}^{z_{k+\halb}} \hspace{-1mm} 
{\bf \tilde f} \! \left({\bf q}_h^-(x_{i+\halb},y,z,t),{\bf q}_h^+(x_{i+\halb},y,z,t)\right) dz \, dy \, dt, 
\end{equation}
\begin{equation}
\label{averG}
{\bf g}_{i,j+\halb,k}=\frac{1}{\Delta t}\frac{1}{\Delta x_i}\frac{1}{\Delta z_k} \hspace{-1mm}  \int \limits_{t^n}^{t^{n+1}} \! \int \limits_{x_{i-\halb}}^{x_{i+\halb}} \! \int \limits_{z_{k-\halb}}^{z_{k+\halb}} \hspace{-1mm} 
{\bf \tilde g} \! \left({\bf q}_h^-(x,y_{j+\halb},z,t),{\bf q}_h^+(x,y_{j+\halb},z,t)\right) dz\,dx\,dt, \\
\end{equation}
\begin{equation}
\label{averH}
{\bf h}_{ij,k+\halb}=\frac{1}{\Delta t}\frac{1}{\Delta x_i}\frac{1}{\Delta y_j} \hspace{-1mm}  \int \limits_{t^n}^{t^{n+1}} \! \int \limits_{x_{i-\halb}}^{x_{i+\halb}} \! \int \limits_{y_{j-\halb}}^{y_{j+\halb}} \hspace{-1mm}  
{\bf \tilde h}\! \left({\bf q}_h^-(x,y,z_{k+\halb},t),{\bf q}_h^+(x,y,z_{k+\halb},t)\right) dy\,dx\,dt  
\end{equation}
are the average fluxes along each Cartesian direction. Furthermore we have defined the space-time average of the smooth part of the non-conservative product as 
\begin{equation}
	{\bf{\bar{P}}}_{ijk}=\frac{1}{\Delta t}\frac{1}{\Delta x_i}\frac{1}{\Delta y_j}\frac{1}{\Delta z_k}\int \limits_{t^n}^{t^{n+1}}\int \limits_{x_{i-\halb}}^{x_{i+\halb}}\int \limits_{y_{j-\halb}}^{y_{j+\halb}}\int \limits_{z_{k-\halb}}^{z_{k+\halb}}{\bf B}({\bf q}_h)\cdot\nabla {\bf q}_h \,\, dz\,\,dy\,\,dx\,\,dt
\end{equation}
and the space-time averaged source term 
\begin{equation}
\label{source:S}
{\bf \bar{S}}_{ijk}=\frac{1}{\Delta t}\frac{1}{\Delta x_i}\frac{1}{\Delta y_j}\frac{1}{\Delta z_k}\int \limits_{t^n}^{t^{n+1}}\int \limits_{x_{i-\halb}}^{x_{i+\halb}}\int \limits_{y_{j-\halb}}^{y_{j+\halb}}\int \limits_{z_{k-\halb}}^{z_{k+\halb}}{ \bf S}\left(\mathbf{q}_h(x,y,z,t)\right) dz\,dy\,dx\,dt\,.
\end{equation}
The terms ${\bf q}_h$ in Eq.~(\ref{averF})--(\ref{source:S}) are piecewise space-time polynomials of degree $M$ and represent the time-evolved reconstruction polynomials.  To obtain the ${\bf q}_h$ first a WENO 
reconstruction polynomial $\w_h$ is obtained from the cell averages ${\bar \u}_{ijk}$ at time $t^n$ (see Sect.\ref{sec:WENO_reconstruction}) and subsequently the time evolution is carried out via a local space-time 
DG predictor as illustrated in Sect.~\ref{sec:LSTDG}. In order to integrate the non-conservative product we use the Dal Maso--Le Floch--Murat theory (see \cite{DLMtheory}) where the non-smooth part of the 
non-conservative term is defined as a Borel measure. In this formulation we therefore also need to account for the jumps of ${\bf q}_h$ at the element boundaries   
\begin{equation}
{{D}}^x_{i+\halb,j,k} \!=\!\! \frac{1}{\Delta t}\frac{1}{\Delta y_j}\frac{1}{\Delta z_k} \hspace{-1mm}  \int \limits_{t^n}^{t^{n+1}} \! \int \limits_{y_{j-\halb}}^{y_{j+\halb}} \! \int \limits_{z_{k-\halb}}^{z_{k+\halb}} \hspace{-2mm}  
{\cal{D}}_1  \! \left({\bf q}_h^-(x_{i+\halb},y,z,t),{\bf q}_h^+(x_{i+\halb},y,z,t)\right) dz \, dy \, dt, 
	$$
	$$
{{D}}^y_{i,j+\halb,k} \!=\!\! \frac{1}{\Delta t}\frac{1}{\Delta x_i}\frac{1}{\Delta z_k} \hspace{-1mm}  \int \limits_{t^n}^{t^{n+1}} \! \int \limits_{x_{i-\halb}}^{x_{i+\halb}} \! \int \limits_{z_{k-\halb}}^{z_{k+\halb}} \hspace{-2mm}  
{\cal{D}}_2  \! \left({\bf q}_h^-(x,y_{j+\halb},z,t),{\bf q}_h^+(x,y_{j+\halb},z,t)\right) dz \, dx \, dt, 
  $$
  $$
{{D}}^z_{i,j,k+\halb} \!=\!\! \frac{1}{\Delta t}\frac{1}{\Delta x_i}\frac{1}{\Delta y_j} \hspace{-1mm}  \int \limits_{t^n}^{t^{n+1}} \! \int \limits_{x_{i-\halb}}^{x_{i+\halb}} \! \int \limits_{y_{j-\halb}}^{y_{j+\halb}} \hspace{-2mm}  
{\cal{D}}_3  \! \left({\bf q}_h^-(x,y,z_{k+\halb},t),{\bf q}_h^+(x,y,z_{k+\halb},t)\right) dy \, dx \, dt, 
\end{equation}
using the path integrals
\begin{equation}
  \label{eqn.dm}
  {\cal{D}}_m({\bf q}_h^-,{\bf q}_h^+) = \int_0^1{{\bf B}_m \left(\Psi({\bf q}_h^-,{\bf q}_h^+,s)\right)\frac{\partial\Psi}{\partial s}ds}, 
\end{equation} 
where $\Psi(s)$ is a path joining the left and right boundary extrapolated states ${\bf q}_h^-$ and ${\bf q}_h^+$ in state space. 
The simplest option is to use a straight-line segment path 
\begin{equation}	\label{segment}
 \Psi = \Psi({\bf q}_h^-,{\bf q}_h^+,s)
={\bf q}_h^- + s({\bf q}_h^+ - {\bf q}_h^-)\,,\qquad 0\leq s \leq 1. 
\end{equation}
Though simplified, 
this particular choice is useful in many applications since in the case of the shallow water equations it guarantees that the resulting numerical scheme is well-balanced and for the 
Baer-Nunziato model it has been shown to preserve the Abgrall condition \cite{AbgrallQuasiCons,AbgrallKarni2} when used with FORCE and Osher-type Riemann solvers \cite{USFORCE2,OsherNC}.  
This may be no longer the case for other systems of equations, for which eventually more sophisticated paths must be adopted, see e.g. \cite{MuellerWB}. 
With the choice of the path (\ref{segment}), the terms ${\cal{D}}_m$ in (\ref{eqn.dm}) can be computed as   
\begin{equation}
\label{Osher-D}
{\cal{D}}_m({\bf q}_h^-,{\bf q}_h^+) = \left( \int_0^1{ {\bf B}_m \left(\Psi({\bf q}_h^-,{\bf q}_h^+,s)\right) ds} \right) \left( {\bf q}_h^+ - {\bf q}_h^- \right). 
\end{equation}
%
%
%
The practical computation of the integrals (\ref{Osher-D}) is typically  performed through 
a three-point Gauss-Legendre formula \cite{USFORCE2,OsherUniversal,OsherNC}. 
Finally, for the numerical approximation of the fluxes (\ref{averF})--(\ref{averH}) we have 
either adopted a local Lax-Friedrichs flux (Rusanov flux) or a simplified Osher--Solomon flux formula 
proposed in \cite{OsherUniversal,OsherNC,ToroOsher}, 
\begin{eqnarray}
\label{eqn.osher}
{\bf \tilde f}\left( \q_h^-, \q_h^+ \right)
=\frac{1}{2}\left( \textbf{f}(\textbf{q}_h^+)+\textbf{f}(\textbf{q}_h^-) \right) + \frac{1}{2} \left( \int \limits_0^1 |{\bf A}_1(\Psi)|ds\, \right) \left(\textbf{q}_h^+-\textbf{q}_h^-\right),
\end{eqnarray}
with 
\begin{equation}
|{\bf A}_1|={\bf R}|{\bf \Lambda}|{\bf R}^{-1}\,,\qquad  |{\bf \Lambda}|={\rm diag}(|\lambda_1|, |\lambda_2|, \ldots, |\lambda_N|)\,,
\end{equation}
and where the path $\Psi$ in Eq. (\ref{eqn.osher}) is the same segment path adopted in (\ref{segment}) for the computation of the jumps ${\cal{D}}_m$. 
Again, the path integral is evaluated numerically using Gauss-Legendre quadrature rules. An entirely analogous procedure allows the computation of the numerical fluxes $\mathbf{\tilde g}$ 
and $\mathbf{\tilde h}$. Note that in the formulation above the numerical flux (\ref{eqn.osher}) contains in its dissipative term both, the conservative and the non-conservative 
part of the system (\ref{NCsyst}).  

\subsection{WENO reconstruction}
\label{sec:WENO_reconstruction}

In order to compute high order intercell fluxes and to integrate the source terms, it is necessary to carry out a \textit{reconstruction} from the cell averages of the solution available  
at the beginning of each time step. To this extent, we provide the basic information about the WENO implementation that we have adopted, which differs from the standard one by \cite{shu_efficient_weno}. 
Although a genuine multidimensional reconstruction is very natural within finite volume methods and is applicable to grids with arbitrary triangulations, 
we have adopted a cheaper dimension-by-dimension methodology \cite{TitarevToroWENO3D,titarevtoro}, which is also much simpler to implement in the presence of adaptive mesh refinement.

A crucial aspect to consider is the choice of the basis functions to be used in the reconstruction process. Two alternative options are available: a \textit{modal} basis or a \textit{nodal} basis,
both of them rescaled on a reference unit interval, e.g. $I=[0;1]$, through the following changes of coordinates valid for each element $I_{ijk}$
\begin{eqnarray}
\label{x-xi}
\xi    &=& \xi(x,i) = \frac{1}{\Delta x_i} \left( x-x_{i-\frac{1}{2}} \right)\,,\\
\label{y-eta}
\eta   &=& \eta(y,j) = \frac{1}{\Delta y_j} \left( y-y_{j-\frac{1}{2}} \right)\,,\\
\label{z-zeta}
\zeta   &=& \zeta(z,k) = \frac{1}{\Delta z_k} \left( z-z_{k-\frac{1}{2}} \right)\,.
\end{eqnarray} 

The modal basis is formed by a set of $M+1$ linearly independent polynomials, typically the Legendre polynomials, having degree from zero to $M$.
The nodal basis is instead formed by a set of $M+1$ linearly independent polynomials, $\{\psi_l\}_{l=1}^{M+1}$, all of degree $M$, 
which are effectively the Lagrange polynomials interpolating a set of $M+1$ nodal points, $\{x_k\}_{k=1}^{M+1}$, in such a way that
\begin{equation}
\psi_l(x_k)=\delta_{lk}\hspace{1cm}l,k=1,2,\ldots, M+1\,.
\end{equation}
Numerical experiments carried out by \cite{HidalgoDumbser} have shown that it is more efficient to use the nodal basis instead of the modal one, especially if the Gauss-Legendre nodes are used.
The reconstruction is performed for each cell $I_{ijk}$ on a reconstruction stencil that, for each Cartesian direction, is given by  
\begin{equation}
\label{eqn.stencildef}  
\mathcal{S}_{ijk}^{s,x} = \bigcup \limits_{e=i-L}^{i+R} {I_{ejk}}, \quad 
\mathcal{S}_{ijk}^{s,y} = \bigcup \limits_{e=j-L}^{j+R} {I_{iek}}, \quad 
\mathcal{S}_{ijk}^{s,z} = \bigcup \limits_{e=k-L}^{k+R} {I_{ije}},  
\end{equation}
where $L$ and $R$, which depend both on the order and on the specific stencil considered, are the spatial extension of the stencil to the left and to the right, respectively. 
Odd order schemes (even polynomial degrees $M$) always adopt three stencils, one central 
stencil ($s=1$, $L=R=M/2$), one fully left--sided stencil ($s=2$, $L=M$, $R=0$) and  one fully  right--sided stencil ($s=3$, $L=0$, $R=M$). 
Even order schemes (odd polynomial degree $M$) always adopt four stencils, two of which 
are central ($s=0$, $L=$floor$(M/2)+1$, $R=$floor$(M/2)$) and ($s=1$, $L=$floor$(M/2)$,  $R=$floor$(M/2)+1$), while the remaining two are again given by the fully left--sided and by the 
fully right--sided stencil, as defined before. The total amount of cells of each stencil is the same as that of the order of the scheme, namely $M+1$. 
%
\subsubsection{ Reconstruction along $x$}
The reconstruction is first performed in the $x$ direction, by writing the reconstructed 
polynomial in terms of the nodal basis $\psi_l(\xi)$ 
\begin{equation}
\label{eqn.recpolydef.x} 
 \w^{s,x}_h(x,t^n) = \sum \limits_{p=0}^M \psi_p(\xi) \hat \w^{n,s}_{ijk,p} := \psi_p(\xi) \hat \w^{n,s}_{ijk,p}\,,
\end{equation}
where we have used the Einstein summation convention, implying summation over indices appearing twice.
Integral conservation on all elements of the stencil then yields the linear equation system from which the 
unknown reconstruction coefficients  $\hat \w^{n,s}_{ijk,p}$ can be determined: 
\begin{equation}
 \frac{1}{\Delta x_e} \int _{x_{e-\halb}}^{x_{e+\halb}} \psi_p(\xi(x)) \hat \w^{n,s}_{ijk,p} \, dx = {\bf \bar u}^n_{ejk}, \qquad \forall {I}_{ejk} \in \mathcal{S}_{ijk}^{s,x}. 
 \label{eqn.rec.x} 
\end{equation}
Once the reconstruction has been performed for each of the stencils relative to the element $I_{ijk}$, we finally construct a data-dependent nonlinear combination of the 
polynomials obtained for each stencil, i.e. 
\begin{equation}
\label{eqn.weno} 
 \w_h^x(x,t^n) = \psi_p(\xi) \hat \w^{n}_{ijk,p}, \quad \textnormal{ with } \quad  
 \hat \w^{n}_{ijk,p} = \sum_{s=1}^{N_s} \omega_s \hat \w^{n,s}_{ijk,p},   
\end{equation}   
where, as specified at the beginning of  Sect.~\ref{sec:WENO_reconstruction}, 
the number of stencils is ${N_s}=3$ or ${N_s}=4$, for even or odd $M$, respectively, while
the nonlinear weights are given by the relations
\begin{equation}
\omega_s = \frac{\tilde{\omega}_s}{\sum_q \tilde{\omega}_q}\,,  \qquad
\tilde{\omega}_s = \frac{\lambda_s}{\left(\sigma_s + \epsilon \right)^r}\,. 
\end{equation} 
The oscillation indicator $\sigma_s$ is
\begin{equation}
\sigma_s = \Sigma_{lm} \hat \w^{n,s}_l \hat \w^{n,s}_m\,,
\end{equation}
and it requires the computation of the oscillation indicator matrix (see \cite{DumbserEnauxToro})
\begin{equation}
\Sigma_{lm} = \sum \limits_{\alpha=1}^M \int \limits_0^1 \frac{\partial^\alpha \psi_l(\xi)}{\partial \xi^\alpha} \cdot \frac{\partial^\alpha \psi_m(\xi)}{\partial \xi^\alpha} d\xi\,. 
\end{equation}
Unlike the original pointwise WENO of \cite{shu_efficient_weno}, which is $2M+1$ order accurate in smooth regions of the solution,
our $M+1$ order accurate implementation of WENO allows for a pragmatic choice of the coefficients $\lambda_s$. In particular,
we select $\lambda_s=1$ for the one--sided stencils and $\lambda_s=10^5$ for the central stencils. Moreover, we use  $\epsilon=10^{-14}$ and $r=8$.

\subsubsection{ Reconstruction along $y$ and $z$}
%
Because the resulting reconstructed polynomial $\w_h^x(x,t^n)$ is only a polynomial in $x$ direction, but still an average in the $y$ and $z$ directions, the reconstruction 
algorithm described above must be applied again along the remaining two directions. In practice, the  
steps from (\ref{eqn.recpolydef.x}) to (\ref{eqn.weno}) are repeated and details about this procedure can be found in \cite{Dumbser2012a}.

\subsection{The local space-time Galerkin predictor}
\label{sec:LSTDG}
The objective of the local space-time Galerkin predictor method is to provide the time evolution, locally for each element, of the reconstructed polynomials $\w_h({\bf x},t^n)$ 
obtained through the WENO reconstruction described before. 
However, unlike the original ADER approach of Titarev and Toro, which was based on a Taylor expansion in time and required a repeated use of the governing 
PDE system in order to substitute time derivatives with space derivatives, the new method relies on a \textit{weak} integral formulation of  
the governing PDE in space--time using an element--local space--time Galerkin method. As a result, all that is required is a point--wise evaluation of fluxes 
and source terms. 
The result of the local space--time Galerkin predictor are the high order space--time polynomials $\q_h$ that are needed for the evaluation of 
the numerical fluxes, source terms and nonconservative jump terms in the scheme (\ref{FVformula}). In the following we briefly illustrate the method, referring to  
\cite{DumbserEnauxToro,Dumbser2008,DumbserZanotti,USFORCE2,HidalgoDumbser,Dumbser2012a} for more details.
 
After introducing the reference time coordinate $\tau=(t-t^n)/\Delta t$, we write the system (\ref{NCsyst}) in terms of the reference coordinates $\tau$ and $\vec \xi = (\xi,\eta,\zeta)$, 
see the definitions (\ref{x-xi})--(\ref{z-zeta}), to get 
\begin{equation}
\label{NCsyst_ref}
\frac{\partial{\bf u}}{\partial \tau} + \frac{\partial \mathbf{f}^\ast}{\partial \xi} + \frac{\partial \mathbf{g}^\ast}{\partial \eta} + \frac{\partial \mathbf{h}^\ast}{\partial \zeta}
+ \mathbf{B}_1^\ast \frac{\partial{\bf u}}{\partial \xi} + \mathbf{B}_2^\ast \frac{\partial{\bf u}}{\partial \eta} + \mathbf{B}_3^\ast \frac{\partial{\bf u}}{\partial \zeta}
 ={\bf S}^\ast \,,
\end{equation}
with
\begin{equation}
{\bf f}^\ast= \frac{\Delta t}{\Delta x_i} \, {\bf f}, \quad 
{\bf g}^\ast= \frac{\Delta t}{\Delta y_j} \, {\bf g}, \quad 
{\bf h}^\ast= \frac{\Delta t}{\Delta z_k} \, {\bf h}, \quad 
\end{equation}
and
\begin{equation}
{\bf B}_1^\ast= \frac{\Delta t}{\Delta x_i} \, {\bf B}_1, \quad 
{\bf B}_2^\ast= \frac{\Delta t}{\Delta y_j} \, {\bf B}_2, \quad 
{\bf B}_3^\ast= \frac{\Delta t}{\Delta z_k} \, {\bf B}_3, \quad 
{\bf S}^\ast= \Delta t {\bf S}. 
\end{equation}

%
%
%
%
%
To obtain a local space-time discontinuous Galerkin approach, we multiply expression (\ref{NCsyst_ref}) by piecewise space-time 
polynomials $\theta_\mathfrak{q}(\xi,\eta,\zeta,\tau)$ of degree $M$, which are given by a tensor--product of the basis functions 
$\psi_l$ adopted in the reconstruction procedure. Here, we use the multi--index $\mathfrak{q}=(p,q,r,s)$. Integration over the 
space-time reference control volume then yields 
\begin{eqnarray}
&& \int \limits_{0}^{1} \int \limits_{0}^{1}  \int \limits_{0}^{1}   \int \limits_{0}^{1}   
\theta_\mathfrak{q} \left( 
  \frac{\partial{\mathbf{q}_h}}{\partial \tau} + \frac{\partial \mathbf{f}^\ast}{\partial \xi} + \frac{\partial \mathbf{g}^\ast}{\partial \eta} + \frac{\partial \mathbf{h}^\ast}{\partial \zeta}   
   \right) d\xi d\eta d\zeta d\tau = \nonumber\\
&&  \int \limits_{0}^{1} \int \limits_{0}^{1}  \int \limits_{0}^{1}   \int \limits_{0}^{1}   
\theta_\mathfrak{q} \left( {\bf S}^\ast - \mathbf{B}_1^\ast \frac{\partial \mathbf{q}_h}{\partial \xi} - \mathbf{B}_2^\ast \frac{\partial \mathbf{q}_h}{\partial \eta} - \mathbf{B}_3^\ast \frac{\partial \mathbf{q}_h}{\partial \zeta} \right) d\xi d\eta d\zeta d\tau, 
\label{eqn.pde.weak1} 
\end{eqnarray}
The discrete space--time solution of Eq. (\ref{eqn.pde.weak1}) will be henceforth denoted by
$\mathbf{q}_h$, and is expanded over the same space-time basis of polynomials as
\begin{equation}
 \mathbf{q}_h(\vec{\xi},\tau) = \theta_\mathfrak{p}(\vec \xi,\tau ) \hat \q_\mathfrak{p},   
 \label{eqn.st.q} 
\end{equation}
where $\hat \q_\mathfrak{p}$ are the unknown \textit{nodal} degrees of freedom. A similar \textit{nodal} representation is provided for the remaining terms entering Eq. (\ref{eqn.pde.weak1}). 
Integration by parts in time yields 
\begin{eqnarray}
 && \int \limits_{0}^{1} \int \limits_{0}^{1}  \int \limits_{0}^{1} \theta_\mathfrak{q}({\vec \xi},1) \theta_\mathfrak{p}({\vec\xi},1) \hat \q_\mathfrak{p}
  \, d\xi d\eta d\zeta 
 - \int \limits_{0}^{1} \int \limits_{0}^{1}  \int \limits_{0}^{1}   \int \limits_{0}^{1} \left(\frac{\partial}{\partial \tau} \theta_\mathfrak{q} \right) \theta_\mathfrak{p} \hat \q_\mathfrak{p}  
 \, d\xi d\eta d\zeta d\tau  = 
   \nonumber \\ 
&&  \phantom{-} \int \limits_{0}^{1} \int \limits_{0}^{1}  \int \limits_{0}^{1} \theta_\mathfrak{q}({\vec \xi},0) \w_h({\vec \xi},t^n) \, d\xi d\eta d\zeta\, \nonumber \\ 
&&  - \int \limits_{0}^{1} \int \limits_{0}^{1}  \int \limits_{0}^{1}   \int \limits_{0}^{1}   
\theta_\mathfrak{q} \left( 
    \frac{\partial}{\partial \xi} \mathbf{f}^\ast(\mathbf{q}_h) + \frac{\partial }{\partial \eta}  \mathbf{g}^\ast(\mathbf{q}_h) + \frac{\partial }{\partial \zeta} \mathbf{h}^\ast  (\mathbf{q}_h)
   \right) d\xi d\eta d\zeta d\tau  \nonumber\\
&&  + \int \limits_{0}^{1} \int \limits_{0}^{1}  \int \limits_{0}^{1}   \int \limits_{0}^{1}   
\theta_\mathfrak{q} \left( {\bf S}^\ast - \mathbf{B}_1^\ast \frac{\partial \mathbf{q}_h}{\partial \xi} - \mathbf{B}_2^\ast \frac{\partial \mathbf{q}_h}{\partial \eta} - \mathbf{B}_3^\ast \frac{\partial \mathbf{q}_h}{\partial \zeta} \right) d\xi d\eta d\zeta d\tau, 
\label{eqn.pde.weak3} 
\end{eqnarray}
Eq. (\ref{eqn.pde.weak3}) should be regarded as a nonlinear algebraic equation, to be solved locally for each element of the computational grid in the unknowns 
$\hat \q_\mathfrak{p}$. Additional details about the local space-time Galerkin predictor for the specific case of non-conservative systems can be found for example in \cite{USFORCE2}.
We emphasize that the choice of a nodal basis based on Gauss-Legendre nodes allows a \textit{dimension-by-dimension} evaluation of the terms appearing in Eq. (\ref{eqn.pde.weak3}). 

\section{Adaptive Mesh Refinement}
\label{sec:AMR}

A detailed illustration of the AMR implementation within our ADER-WENO approach has been presented
in \cite{Dumbser2012a}. The description we provide here is self-contained, but focused on the essential aspects only.
Unlike the original patch-based block-structured approach by Berger \& Oliger~\cite{Berger-Oliger1984,Berger-Jameson1985,Berger-Colella1989}, we 
refine individual Cartesian cells, which are treated as elements of a tree data structure, like in \cite{Khokhlov1998}.  This choice is particularly 
suited to our element-local space-time DG predictor, which does not need any exchange of information through neighbor elements, and can therefore be implemented
with no modifications even if two adjacent cells belong to two different levels of grid refinement. 

\subsection{AMR implementation}
Each level of refinement is indicated with $\ell$, ranging from the coarsest level $\ell=0$
to the maximally refined level $\ell=\ell_{\rm max}$, beyond which no further refinement is possible.
In addition, we use ${\mathcal L}_{\ell}$ to denote the union of all elements 
up to level $\ell$. 
Any cell ${\mathcal C}_m$, at any level of refinement, is identified with a unique positive
integer number $m$, with $1\leq m \leq N_{{\rm Cells}}$, where 
$N_{{\rm Cells}}$ is the (time--dependent) total number of cells at any given time.
When a cell  ${\mathcal C}_m$ at level $\ell$ is refined, we refer to it as a 
\textit{mother} or parent cell, and the cells on the next refinement level $\ell+1$ contained in it are called \textit{children} cells. 
Moreover, the \textit{Neumann neighbors} $\mathcal{N}_m$ of a cell ${\mathcal C}_m$ are the neighbor cells that share a 
face with ${\mathcal C}_m$. Each cell has $2d$ Neumann neighbors in $d$ space  dimensions, except of course the case of
the cells at the boundaries of the computational domain. 
On the other hand, the \textit{Voronoi neighbors} $\mathcal{V}_m$ of a cell ${\mathcal C}_m$ are those cells which 
share  at least one node with ${\mathcal C}_m$, and each cell has $3^d-1$ Voronoi neighbors in $d$ space dimensions.    

\begin{figure}
\begin{center}
\includegraphics[scale=.4,angle=-90]{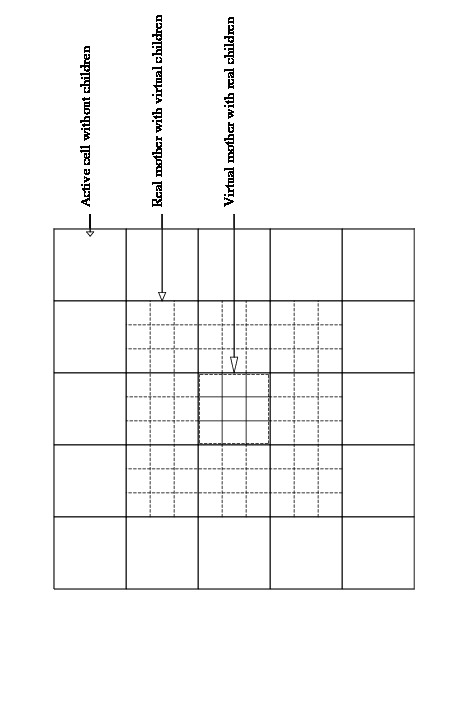}
\caption{Sketch of a cell-by-cell refinement of the central cell ${\mathcal C}_m$ on level $\ell$. The children of ${\mathcal C}_m$ are real cells on level $\ell+1$, 
surrounded by \textit{virtual} children of the Voronoi neighbors $\mathcal{V}_m$ of cell ${\mathcal C}_m$. The virtual children are on the same level $\ell+1$ and are
needed for reconstruction.}
\label{fig_AMRCells}
\end{center}
\end{figure}

Since the high order finite volume schemes used in this paper need information from neighbors to carry out the WENO reconstruction (even more than just 
the Voronoi neighbors $\mathcal{V}_m$), each refined cell on a level $\ell+1$ must be surrounded by a layer of either real or \textit{virtual} cells on 
the same level. The layer thickness must be greater or equal to the size of the reconstruction stencil. Likewise, a mother cell on level $\ell$ that is refined 
continues to exist as a virtual mother cell since it may be surrounded by non-refined cells on the same level which need its information for reconstruction. 
A schematic representation of this mechanism involving one level of refinement is reported in figure~\ref{fig_AMRCells}. There, the central cell of level 
$\ell$ is refined, hence it becomes virtual and has real children on level $\ell+1$, while the surrounding cells are virtually refined in order to allow
the real cells on level $\ell+1$ to perform the WENO reconstruction. 

Having introduced this terminology, it
is convenient to list schematically the rules that are adopted in our cell-by-cell AMR implementation.
\begin{itemize}
\item 
Any AMR scheme requires a criterion for deciding whether a given cell ${\mathcal C}_m$ needs refinement or recoarsening.
We have adopted the same strategy described in 
\cite{Lohner1987}, which is based on the  calculation of 
a second derivative error.
A cell ${\mathcal C}_m$ is marked for refinement  if $\chi_m>\chi_{\rm ref}$, while it is marked for 
recoarsening if $\chi_m<\chi_{\rm rec}$, where
\begin{equation}
\chi_m=\sqrt{\frac{\sum_{k,l} (\partial^2 \Phi/\partial x_k \partial x_l)^2 }{\sum_{k,l}[(|\partial \Phi/\partial x_k|_{i+1}+|\partial \Phi/\partial x_k|_i)/\Delta x_l+\varepsilon|(\partial^2 /\partial x_k \partial x_l )||\Phi|]^2} }\,.
\label{eqn.indicator}
\end{equation}
The summation $\sum_{k,l}$ is taken over the number of space dimension of the problem in order to include the cross term derivatives, whereas $\Phi=\Phi(\u)$  
is a generic indicator function of the conservative variables $\u$. In most cases 
we have adopted $\chi_{\rm ref}$ in the range $\sim[0.2,0.25]$ and $\chi_{\rm rec}$ in the 
range $\sim[0.05,0.15]$. Finally, the parameter $\varepsilon$ acts as a filter which prevents refinement in regions of small  ripples and is given the value $\varepsilon = 0.01$.

\item 
Whenever a {\em mother} cell of the level $\ell$ is refined, it generates $\mathfrak{r}^d$ {\em children} cells,
such that 
\begin{equation}
\label{refine-factor}
\Delta x_{\ell} = \mathfrak{r} \Delta x_{\ell+1}\, \quad \Delta y_{\ell} = \mathfrak{r} \Delta y_{\ell+1} \, \quad 
\Delta z_{\ell} = \mathfrak{r} \Delta z_{\ell+1}. 
\end{equation}
In addition, and as commented below,
the time steps are also chosen \textit{locally} on each level so that 
\begin{equation}
\Delta t_{\ell} = \mathfrak{r} \Delta t_{\ell+1}.  
\end{equation}
Due to the high order WENO reconstruction, the refinement factor $\mathfrak{r}$ must satisfy $\mathfrak{r} \geq M$.  
\item 
At any level of refinement,
each cell ${\mathcal C}_m$ has one among three possible 
\textit{status} flags, which we denote by $\sigma$. The first possibility is that ${\mathcal C}_m$  
is an \textit{active cell} ($\sigma=0$), in which case it is updated 
through the finite volume scheme described in the previous Section \ref{sec:Numscheme}. 
The second possibility is that ${\mathcal C}_m$ is
a \textit{virtual child cell} ($\sigma=1$) and is updated by projection of the mother's high order space--time polynomial.
In practice, the virtual children  receive their values
from the active mother via standard $L_2$ projection. For this purpose, the space--time polynomials $\q_h$ can be 
conveniently evaluated at any time. This operation is needed for performing the reconstruction on 
the finer grid level at intermediate times. The projection operator for a cell $\mathcal{C}_m$ on
level $\ell$ is simply given by evaluating the space--time polynomial $\q_h$ of its \textit{mother}  
at any given time $t^n_\ell$ as follows: 
\begin{equation}
  \bar \u_m(t^n_\ell) = \frac{1}{\Delta x_\ell} \frac{1}{\Delta y_\ell}  \frac{1}{\Delta z_\ell}  \int \limits_{\mathcal{C}_m} \q_h(\mathbf{x},t^n_\ell) d \mathbf{x}.  
\end{equation} 
Finally, ${\mathcal C}_m$ can be a \textit{virtual mother cell} ($\sigma=-1$), updated by recursively  averaging over all children from higher refinement levels.  Namely, the virtual mother cell obtains 
its cell average by averaging recursively over the cell averages of all its children, namely 
including the possible children of their children. 
If we denote the set of children of a cell $\mathcal{C}_m$ by $\mathcal{B}_m$,
then the averaging operator is given by 
\begin{equation}
  \bar \u_m = \frac{1}{\mathfrak{r}^d} \sum \limits_{\mathcal{C}_k \in \mathcal{B}_m} \bar \u_k. 
\label{eqn.average} 
\end{equation}  
%
\item 
Only real cells ($\sigma=0$) can be refined. Therefore, if a virtual cell needs to be refined it must be first activated.
\item 
The levels of refinement of two cells that are Voronoi neighbors of each other can only 
differ by at most unity. Moreover, every cell has  Voronoi neighbors, which can be either active or virtual, at the same level of refinement. 
%
\end{itemize}

\subsection{AMR local time stepping}
%
%
\begin{figure}
\begin{center}
\includegraphics[scale=.5]{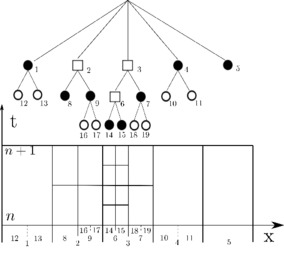}
\includegraphics[scale=.5]{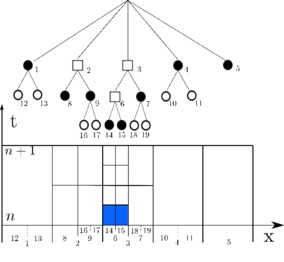}
\includegraphics[scale=.5]{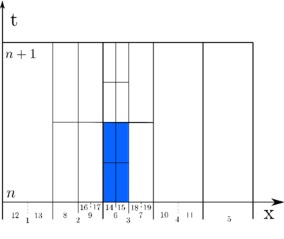}
\includegraphics[scale=.5]{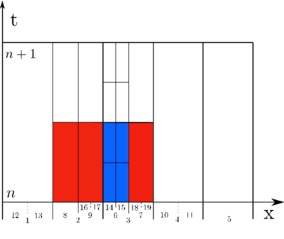}
\includegraphics[scale=.5]{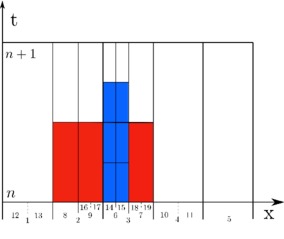}
\includegraphics[scale=.5]{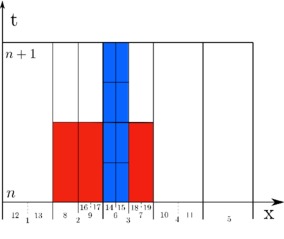}
\includegraphics[scale=.5]{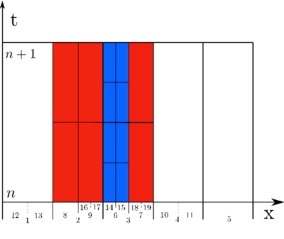}
\includegraphics[scale=.5]{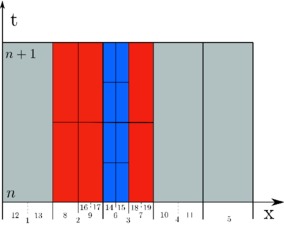}
\caption{Example of a local time stepping algorithm involving two levels of refinement with $\mathfrak{r}=2$. Legend: $\CIRCLE$ (regular active cell), $\Circle$ (virtual refined cell), $\square$ (virtual coarse cell).}
\label{fig_LocalTimeStepping}
\end{center}
\end{figure}
As anticipated above,
each of the refinement levels is advanced in time with its own \textit{local} time-step, i.e. 
$\Delta t_\ell = \mathfrak{r} \Delta t_{\ell+1}$. The use of
time steps that are integer multiples of each other among subsequent levels is very convenient, and indeed very natural within
AMR. However, alternative local time stepping schemes are also possible  
(see \cite{dumbserkaeser06d,TaubeMaxwell,stedg1}), where a different local time step is 
allowed for each element. 
After denoting by $t^n_\ell$ and $t^{n+1}_\ell$ the current and future times of the level $\ell$, 
the first level to be updated is the one with the largest value of $\ell$ satisfying the  
\textit{update criterion}\footnote{We define $t^{n+1}_{-1}:=t^{n+1}_{0}$, so that also the scheduling of 
level $\ell=0$ is taken into account.}~\cite{dumbserkaeser06d}  
\begin{equation} 
t^{n+1}_{\ell} \leq t^{n+1}_{\ell-1},  \qquad  0 \leq \ell \leq \ell_{\max}\,.
\label{eqn.update.criterion}  
\end{equation} 
In practice, starting from the common initial time $t=0$, the finest level of refinement  $\ell_{\max}$ is evolved first and performs a number of $\mathfrak{r}$ sub-timesteps before 
the next coarser level $\ell_{\max}-1$ performs its first time update. This procedure is then applied recursively and it implies a total amount of $\mathfrak{r}^\ell$ sub-timesteps on 
each level to be performed in order to reach the time $t_0^{n+1}$ of the coarsest level. As example of the application of the local time stepping strategy
to a one dimensional case involving two levels of refinement with $\mathfrak{r}=2$
is reported in Fig.\ref{fig_LocalTimeStepping}.

Thanks to the use of the local space--time predictor, 
which computes the predictor solution $\q_h$ for each element 
valid from time $t^n_\ell$ to time $t^{n+1}_\ell$,
the computation of numerical fluxes between two adjacent cells on different levels of 
refinement is rather straightforward.
Further details about the  actual  implementation 
of the local time stepping procedure and of the AMR parallelization
through the standard Message Passing Interface (MPI)
can be found in \cite{Dumbser2012a}. 
%

\section{The Baer-Nunziato equations}
\label{sec:BN}

A Baer-Nunziato type model for compressible two-phase flow without relaxation terms is given by the following 
system of equations, see \cite{BaerNunziato1986,SaurelAbgrall,AndrianovWarnecke,Schwendeman,DeledicquePapalexandris,MurroneGuillard}:

\begin{equation}\left.
\renewcommand{\arraystretch}{1.5}
\begin{array}{l}
\label{ec.BN}


	\frac{\partial}{\partial t}\left(\ar\right)+\nabla\cdot\left(\ar\ub\right)=0,
	
	\\


\frac{\partial}{\partial t}\left(\ar\ub\right)
+\nabla\cdot\left(\phi_1\rho_1\textbf{v}_1 \textbf{v}_1\right)+\nabla\phi_1p_1 
= p_I\nabla\phi_1, 

\\


\frac{\partial}{\partial t}\left(\phi_1\rho_1E_1\right)
+\nabla\cdot\left(\left(\phi_1\rho_1E_1+\phi_1p_1\right)\textbf{v}_1\right) = 
-p_I\partial_t\phi_1, 

\\


\frac{\partial}{\partial t}\left(\phi_2\rho_2\right)+\nabla\cdot\left(\phi_2\rho_2\textbf{v}_{2}\right)=0,

\\ 


\frac{\partial}{\partial t}\left(\phi_2\rho_2\textbf{v}_{2}\right)
+\nabla\cdot\left(\phi_2\rho_2\textbf{v}_2 \textbf{v}_{2}\right)+\nabla\phi_2p_2=p_I\nabla\phi_2,   

\\


\frac{\partial}{\partial t}\left(\phi_2\rho_2E_2\right)
+\nabla\cdot\left(\left(\phi_2\rho_2E_2+\phi_2p_2\right)\textbf{v}_{2}\right)=p_I\partial_t\phi_1, 

\\


\frac{\partial}{\partial t}\phi_1+\textbf{v}_{I}\nabla\phi_1 = 0. 
\end{array}\right\}
\end{equation}

The system is closed by the stiffened equation of state (EOS): 
\begin{equation}
\label{eqn.eos} 
   e_k = \frac{p_k + \gamma_k \pi_k}{\rho_k (\gamma_k -1 )}. 
\end{equation}
Here, $\phi_k$ denotes the volume fraction of phase $k$, $\rho_k$ is the 
density, $\mathbf{v}_k$ is the velocity vector, 
$E_k = e_k + \halb \mathbf{v}_k^2$ and $e_k$ are the phase specific total 
and internal energies, respectively.
In the literature, one of the phases is often also called the \textit{solid} 
phase and the other one the \textit{gas} phase. Defining arbitrarily the 
first phase as the solid phase in the rest of the paper we will therefore use 
the subscripts $1$ and $s$ as well as $2$ and $g$ as synonyms. 
For the interface velocity $\mathbf{u}_I$ and pressure  $p_I$ we choose 
$\mathbf{v}_I = \mathbf{v}_1$ and $p_I = p_2$, according to
 \cite{BaerNunziato1986}, although other choices are also possible 
(see e.g. the paper by Saurel and Abgrall \cite{SaurelAbgrall}). 
The state vector $\u$ is 
\begin{equation}
\u=\left(
\ar, \, \ar \mathbf{v}_1, \, \ar E_1, \, 
\arr, \, \arr \mathbf{v}_2, \, \arr E_2, \, \phi_1 
\right).   
\end{equation}
The system (\ref{ec.BN}) can be cast in the form prescribed by (\ref{NCsyst}) by collecting all the non-conservative terms in the matrix 
$\bf B(\u)$, while keeping the conservative part of the system expressed through $\textbf{F}(\u)$. An exhaustive treatment of the mathematical 
properties of the Baer-Nunziato equations can be found in \cite{AndrianovWarnecke,DeledicquePapalexandris,Schwendeman,TokarevaToro}, where also
exact and approximate solutions to the Riemann problem are given. 


%
%

\section{Test Problems}
\label{sec:tests}

In all test problems shown below the indicator function for the refinement and recoarsening criterion has been chosen as 
\begin{equation} 
\Phi = \sqrt{ \left( \frac{\phi_1 \rho_1}{\rho_{1,0}} \right)^2 + \left( \frac{\phi_2 \rho_2}{\rho_{2,0}} \right)^2 }, 
\end{equation} 
with some reference densities $\rho_{1,0}$ and $\rho_{2,0}$, respectively.  
\subsection{Smooth vortex problem} 
The first test that we have considered is given by a stationary and axisymmetric solution of the Baer-Nunziato equations and has first been reported in \cite{USFORCE2}.  
The resulting configuration describes a vortex-type solution with no radial motion. Because of these assumptions, the continuity and the energy equations are 
automatically satisfied. After choosing a simple dependence on the radius $r$ of the pressure and of the volume fraction of the solid phase, namely 
\begin{eqnarray}
\label{sol.BN.pres}
p_k&=&p_{k0}\left(1-\frac{1}{4}e^{\displaystyle{\left(1-r^2/s_k^2\right)}}\right), 
 \qquad (k=1,2)\,,\\
\label{sol.BN.alph}
	\phi_1&=&\frac{1}{3}+\frac{1}{2\sqrt{2\pi}}e^{\displaystyle{-r^2/2}}\,,
\end{eqnarray}
the momentum equations can be easily solved, to provide the velocity field of the vortex as
\footnotesize
\begin{eqnarray}
\label{sol.BN.vel}	
u_1^\theta&=&\displaystyle{\frac{1}{2s_1D}}
\sqrt{rD\left[p_{10}\left(4\sqrt{2\pi}F_1+6H_1-12Gs_1^2+3H_1s_1^2\right)+3p_{20}s_1^2\left(4G-H_2\right)\right]}\,,\\
u_2^\theta&=&\displaystyle{ \frac{r\sqrt{2}}{2\rho_2s_2} } 
\sqrt{\rho_2p_{20}F_2}\; ,
\end{eqnarray}
\normalsize
where
\begin{equation}	
H_k=e^{\displaystyle{-\frac{2r^2+r^2s_k^2-2s_k^2}{2s_k^2}}}, \quad 
F_k=e^{\displaystyle{-\frac{(r-s_k)(r+s_k)}{s_k^2}}}, \quad  (k=1,2),
\end{equation}
and
\begin{equation}
G=e^{{-{r^2}/{2}}}, \qquad 
D=\rho_1\left(2\sqrt{2\pi}+3G\right).	
\end{equation}
In order to make the test problem unsteady, the vortex is then boosted along the diagonal of the computational domain through a Galilean 
transformation of the velocity, with components $\bar{u}=\bar{v}$. In our tests we have chosen the following parameters 
\begin{equation}
\label{eqn.BN.param}
 \rho_1 = 1, \quad \rho_2 = 2, \quad p_{10} = 1, \quad p_{20} = \frac{3}{2}, \quad s_1=\frac{3}{2}, \quad s_{2}=\frac{7}{5}, 
 \quad \bar{u}=\bar{v}=2\,,  
\end{equation}
while the computational domain is $\Omega=[-10;10] \times [-10;10]$  with four  
periodic boundary conditions, in such a way that the 
exact solution of the problem is given by the initial condition after $T=10$. 
In Table~\ref{tab.conv2} we have reported the results of the convergence tests, where we have used the third and fourth order version of the method.
Here $\rho_{1,0}=\rho_{2,0}=1$ have been chosen.  
The convergence rates have been computed with respect to an initially uniform mesh, as proposed by Berger and Oliger in \cite{Berger-Oliger1984}. 
\begin{table}[!t]   
\caption{Numerical convergence results for the vortex test using the one--step ADER-WENO finite volume scheme.
The error norms refer to the variable $\rho_1$ at the final time $T=10$, and have been computed with $\ell_{\rm max}=2$.
The asterisk $^\ast$ refers to a uniform grid. 
}
\begin{center} 
\renewcommand{\arraystretch}{1.0}
\begin{tabular}{lccllc} 
\hline
  $N_G\times N_G$  & $\epsilon_{L_2}$ & $\mathcal{O}(L_2)$ & $N_G\times N_G$ & $\epsilon_{L_2}$ & $\mathcal{O}(L_2)$  \\ 
\hline
                     & & {$\mathcal{O}3$} & & & {$\mathcal{O}4$}  \\
\hline                     
 15$\times$15$^\ast$   & 4.9627E-01 &      & 15$\times$15$^\ast$  & 4.6443E-01 &       \\ 
 30$\times$30          & 2.5428E-02 & 4.29 & 30$\times$30         & 2.3166E-02 & 4.33   \\ 
 45$\times$45          & 1.3665E-02 & 3.27 & 45$\times$45         & 1.0674E-02 & 3.43   \\ 
 60$\times$60          & 7.8621E-03 & 2.99 & 60$\times$60         & 1.0115E-03 & 4.42   \\  
 90$\times$90          & 2.0279E-03 & 3.07 & 75$\times$75         & 5.6484E-04 & 4.17   \\
 120$\times$120        & 9.9613E-04 & 2.99 & 90$\times$90         & 2.9489E-04 & 4.11   \\
\hline 
\end{tabular} 
\end{center}
\label{tab.conv2}
\end{table} 
\begin{figure}
\begin{center}
\includegraphics[angle=0,width=4.3cm]{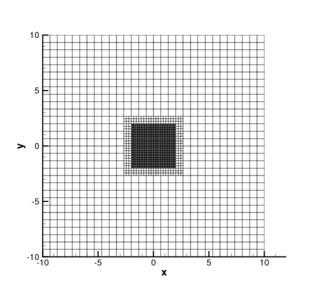}
\includegraphics[angle=0,width=4.3cm]{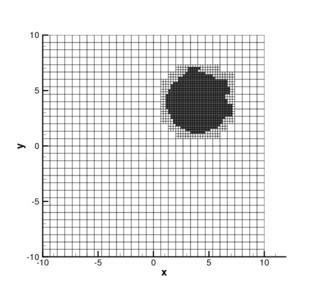}
\includegraphics[angle=0,width=4.3cm]{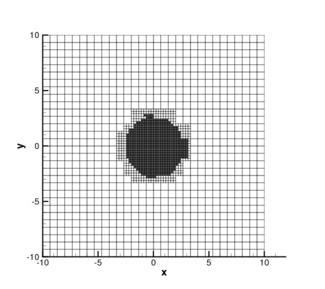}
\caption{  
{
AMR grid of the isentropic-vortex test at the initial time (left panel), at time $t=2.0$ (central panel) and at  
the final time $t=10$ (right panel)}.  
Two levels of refinement have been adopted ($\ell_{\max}=2$), starting from a uniform $45\times45$ grid. 
}
\label{fig_shu_vortex}
\end{center}
\end{figure}
\subsection{Riemann problems}
\label{sec.bnrp} 

The high order space-time adaptive ADER-WENO methods proposed in this paper are particularly
well suited for the accurate resolution of material interfaces in multi-fluid and multi-phase 
flow problems. This claim is validated in the following by solving a set of 1D shock tube problems on 
space-time adaptive Cartesian meshes in 2D. The exact solutions for the Riemann problems solved 
here have been provided in \cite{AndrianovWarnecke,Schwendeman,DeledicquePapalexandris}. 
From these papers a set of six Riemann problems has been chosen, see Table \ref{tab.rpbn.ic}. The 
same set of problems has already been solved with high order unstructured one-step WENO finite volume 
schemes using a centered path-conservative FORCE method in \cite{USFORCE2}. A subset of these Riemann 
problems has been solved again in \cite{OsherNC} using the more accurate path-conservative extension 
of the Osher method, which is also used in this paper. 

The two-dimensional computational domain 
is $\Omega=[-0.5;0.5] \times [0; 1]$, which is discretized at the level 0 grid with only $100 \times 10$ cells.  
A maximum number of two refinement levels ($\ell_{\rm max} = 2$) is chosen, together with a refinement factor 
of $\mathfrak{r}=4$. The discontinuity is located initially at $x=0$ and the final simulation times are given 
in Table \ref{tab.rpbn.ic}. For all test cases a third order ADER-WENO scheme is used with reconstruction in 
characteristic variables. In $x$-direction transmissive boundary conditions are imposed and periodic boundaries 
are applied in $y$-direction. The parameters for the refinement criterion are $\rho_{1,0}=\rho_{2,0}=1$ apart for
RP2 and RP4, where $\rho_{1,0}=1000$ and $\rho_{2,0}=1$.   

\begin{figure}[!htbp]
\begin{center}
\begin{tabular}{cc} 
\includegraphics[width=0.35\textwidth]{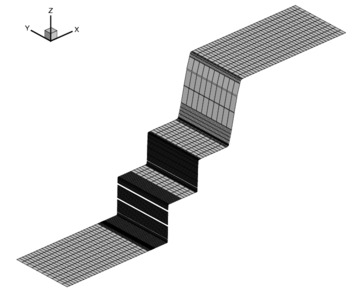}      & 
\includegraphics[width=0.35\textwidth]{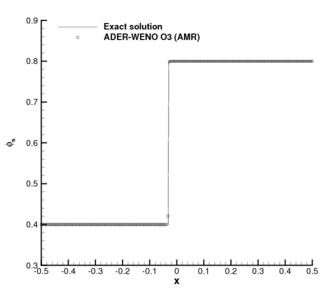}    \\ 
\includegraphics[width=0.35\textwidth]{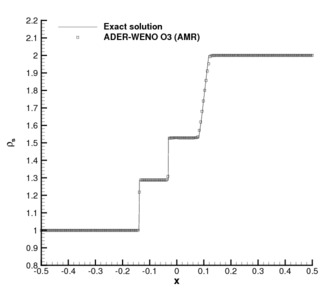}    & 
\includegraphics[width=0.35\textwidth]{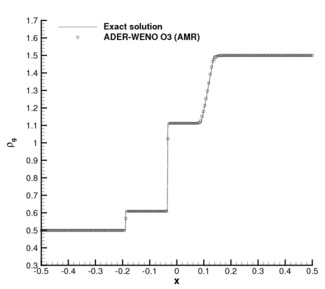}    \\ 
\includegraphics[width=0.35\textwidth]{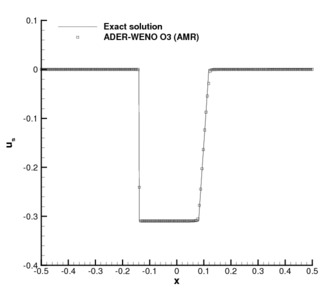}      &  
\includegraphics[width=0.35\textwidth]{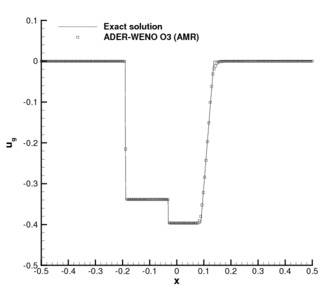}      \\ 
\includegraphics[width=0.35\textwidth]{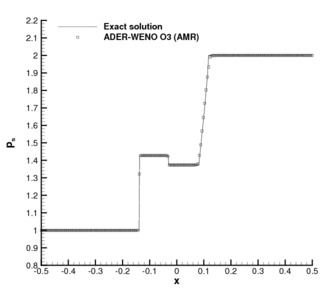}      & 
\includegraphics[width=0.35\textwidth]{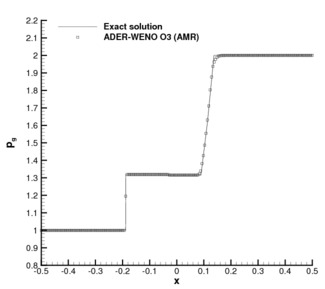}   
\end{tabular}
\caption{Results for Riemann problem RP1.}
\label{fig.bn.rp1}
\end{center}
\end{figure}

\begin{figure}[!htbp]
\begin{center}
\begin{tabular}{cc} 
\includegraphics[width=0.35\textwidth]{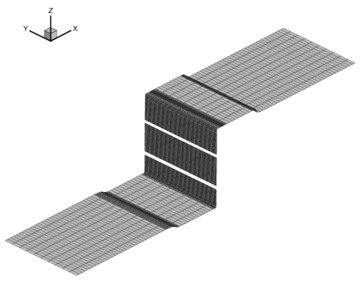}      & 
\includegraphics[width=0.35\textwidth]{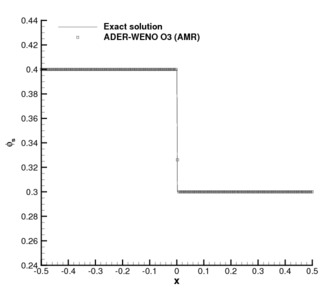}    \\ 
\includegraphics[width=0.35\textwidth]{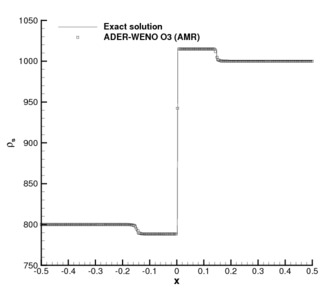}    & 
\includegraphics[width=0.35\textwidth]{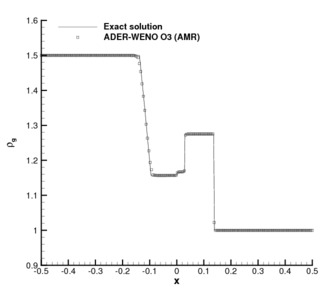}    \\ 
\includegraphics[width=0.35\textwidth]{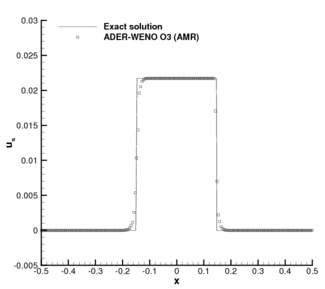}      &  
\includegraphics[width=0.35\textwidth]{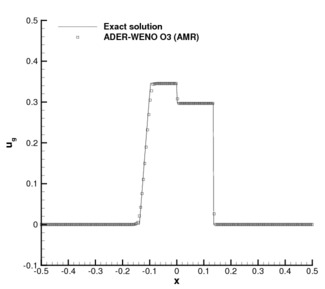}      \\ 
\includegraphics[width=0.35\textwidth]{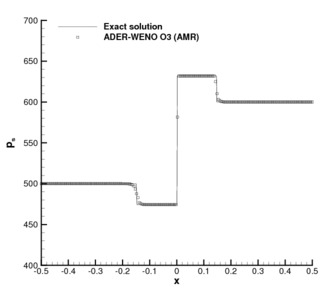}      & 
\includegraphics[width=0.35\textwidth]{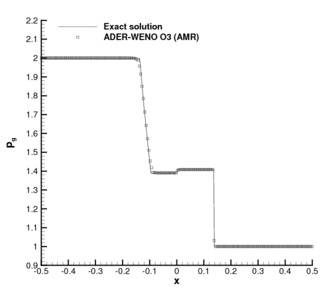}   
\end{tabular}
\caption{Results for Riemann problem RP2.}
\label{fig.bn.rp2}
\end{center}
\end{figure}

\begin{figure}[!htbp]
\begin{center}
\begin{tabular}{cc} 
\includegraphics[width=0.35\textwidth]{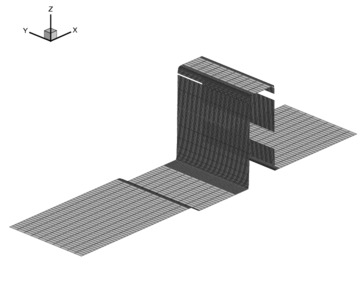}      & 
\includegraphics[width=0.35\textwidth]{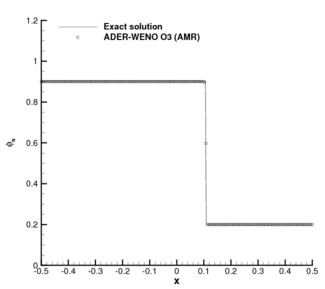}    \\ 
\includegraphics[width=0.35\textwidth]{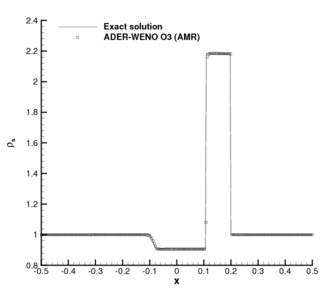}    & 
\includegraphics[width=0.35\textwidth]{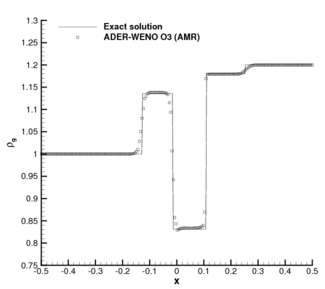}    \\ 
\includegraphics[width=0.35\textwidth]{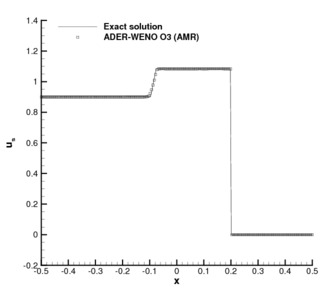}      &  
\includegraphics[width=0.35\textwidth]{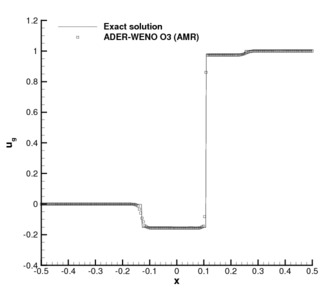}      \\ 
\includegraphics[width=0.35\textwidth]{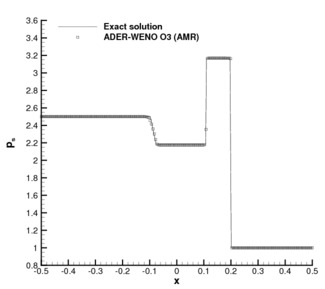}      & 
\includegraphics[width=0.35\textwidth]{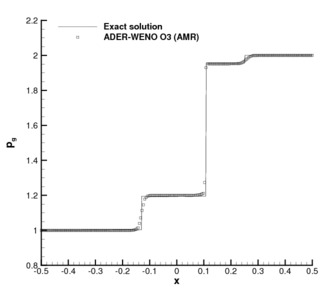}   
\end{tabular}
\caption{Results for Riemann problem RP3.}
\label{fig.bn.rp3}
\end{center}
\end{figure}

\begin{figure}[!htbp]
\begin{center}
\begin{tabular}{cc} 
\includegraphics[width=0.4\textwidth]{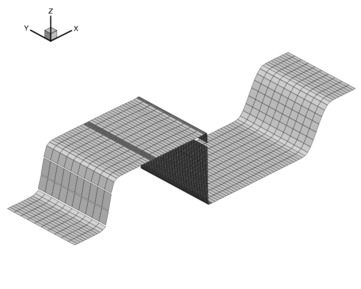}       & 
\includegraphics[width=0.4\textwidth]{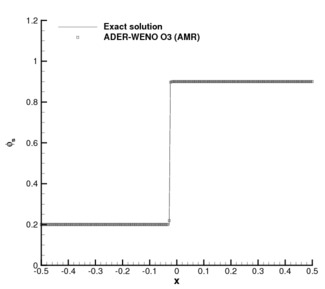}     \\ 
\includegraphics[width=0.4\textwidth]{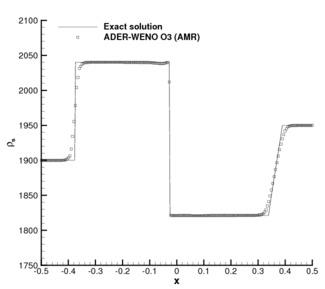}     & 
\includegraphics[width=0.4\textwidth]{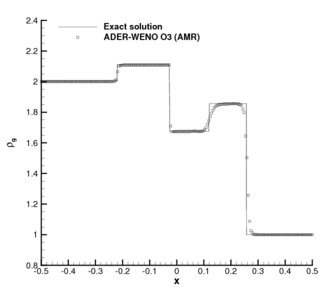}     \\ 
\end{tabular}
\caption{Results for Riemann problem RP4.}
\label{fig.bn.rp4}
\end{center}
\end{figure}

\begin{figure}[!htbp]
\begin{center}
\begin{tabular}{cc} 
\includegraphics[width=0.35\textwidth]{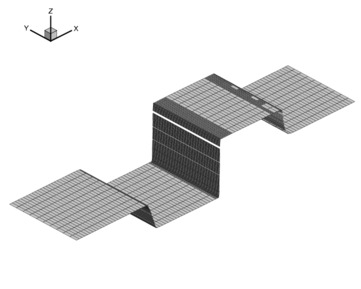}       & 
\includegraphics[width=0.35\textwidth]{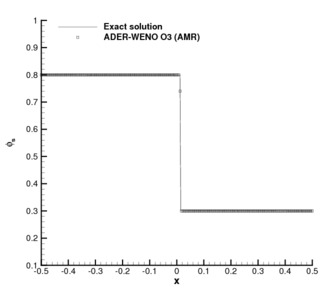}    \\ 
\includegraphics[width=0.35\textwidth]{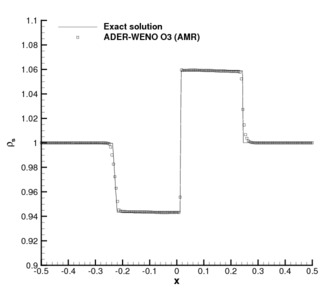}    & 
\includegraphics[width=0.35\textwidth]{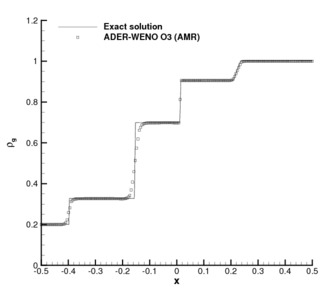}    \\ 
\includegraphics[width=0.35\textwidth]{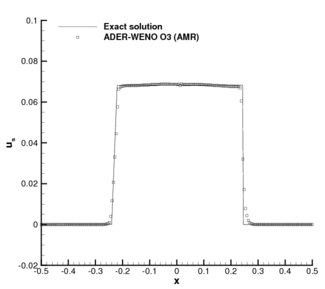}      &  
\includegraphics[width=0.35\textwidth]{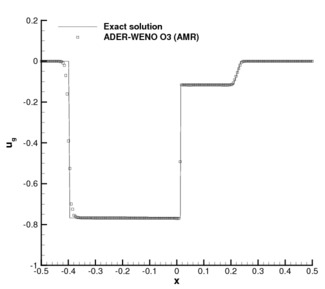}      \\ 
\includegraphics[width=0.35\textwidth]{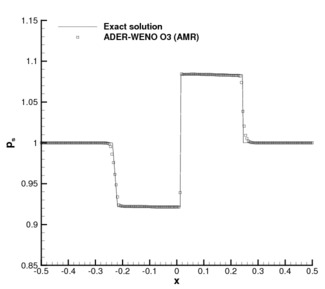}      & 
\includegraphics[width=0.35\textwidth]{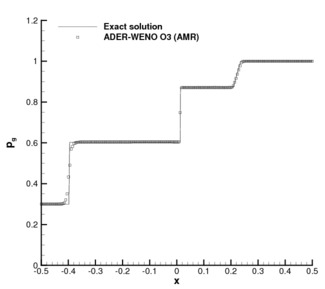}   
\end{tabular}
\caption{Results for Riemann problem RP5.}
\label{fig.bn.rp5}
\end{center}
\end{figure}

\begin{figure}[!htbp]
\begin{center}
\begin{tabular}{cc} 
\includegraphics[width=0.35\textwidth]{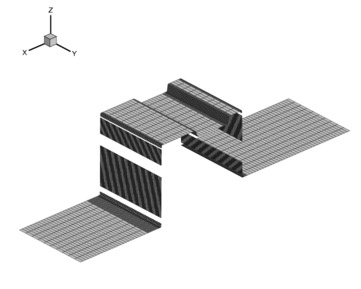}       & 
\includegraphics[width=0.35\textwidth]{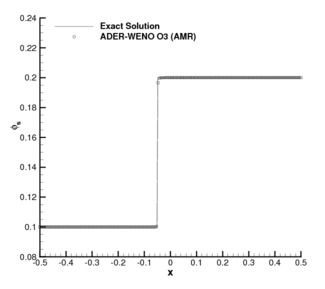}    \\ 
\includegraphics[width=0.35\textwidth]{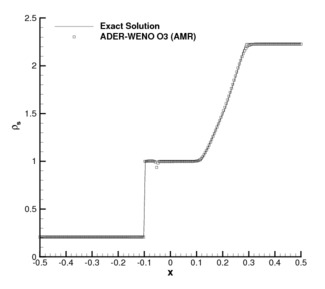}    & 
\includegraphics[width=0.35\textwidth]{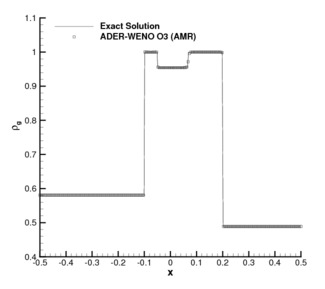}    \\ 
\includegraphics[width=0.35\textwidth]{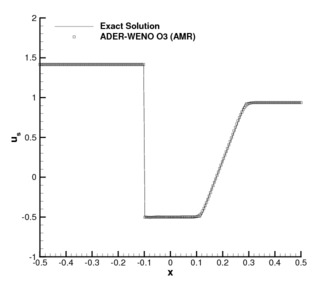}      &  
\includegraphics[width=0.35\textwidth]{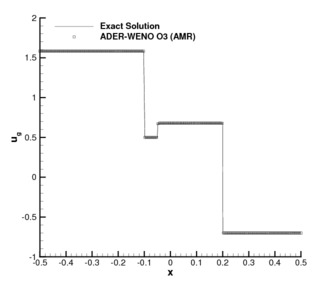}      \\ 
\includegraphics[width=0.35\textwidth]{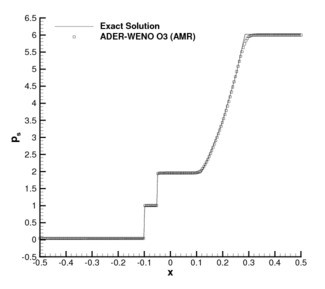}      & 
\includegraphics[width=0.35\textwidth]{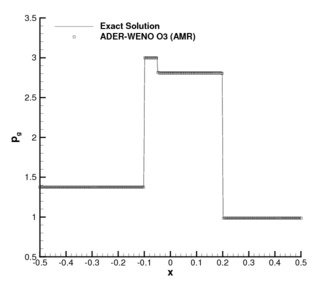}   
\end{tabular}
\caption{Results for Riemann problem RP6.}
\label{fig.bn.rp6}
\end{center}
\end{figure}

The results for the AMR computations are shown in Figs. \ref{fig.bn.rp1} - \ref{fig.bn.rp6}. A sketch of the 
AMR grid at the final simulation time is depicted on the top left of each figure together with a cut through 
the reconstructed numerical solution $\w_h$ on 200 \textit{equidistant} points along the $x$-axis in the remaining 
sub-figures. For RP4 the same quantities as in \cite{Schwendeman} are shown. For all six problems we obtain an 
excellent agreement between the high order AMR computations and the exact reference solutions provided in  
\cite{AndrianovWarnecke,Schwendeman,DeledicquePapalexandris}. The agreement is much better than the one obtained 
in the previous publications of the authors \cite{USFORCE2,OsherNC}. The solid contact is resolved perfectly well in 
all cases.  
Furthermore, no spurious post shock oscillations as reported in \cite{USFORCE2} for RP5 are visible in the present
high order ADER-WENO simulations \textit{with} AMR. Our results clearly confirm that the combination of adaptive mesh
refinement (AMR) with high order WENO finite volume schemes with the little diffusive Osher-type Riemann solver 
\cite{OsherUniversal,OsherNC} is very well suited for the simulation of compressible multi-phase flows, as claimed 
at the beginning of this section. 

\begin{table}[!b]
\caption{Initial states left (L) and right (R) for the Riemann problems solved in 2D. 
Values for $\gamma_i$, $\pi_i$ and the final time $t_e$ are also given.} 
\begin{center}
\renewcommand{\arraystretch}{1.0}
\begin{tabular}{ccccccccc}
\hline
   & $\rho_s$ & $u_s$  & $p_s$ & $\rho_g$ & $u_g$ & $p_g$ & $\phi_s$ & $t_e$  \\
\hline 
\multicolumn{1}{l}{\textbf{RP1 \cite{DeledicquePapalexandris}:} } & 
\multicolumn{8}{c}{ $\gamma_s = 1.4, \quad \pi_s = 0, \quad \gamma_g = 1.4, \quad \pi_g = 0$}  \\
\hline 
L & 1.0    & 0.0   & 1.0  & 0.5 & 0.0   &  1.0 & 0.4 & 0.10 \\
R & 2.0    & 0.0   & 2.0  & 1.5 & 0.0   &  2.0 & 0.8 &      \\
\hline 
\multicolumn{1}{l}{\textbf{RP2 \cite{DeledicquePapalexandris}:}} & 
\multicolumn{8}{c}{ $\gamma_s = 3.0, \quad \pi_s = 100, \quad \gamma_g = 1.4, \quad \pi_g = 0$}  \\
\hline
L & 800.0   & 0.0   & 500.0  & 1.5 & 0.0   & 2.0 & 0.4 & 0.10  \\
R & 1000.0  & 0.0   & 600.0  & 1.0 & 0.0   & 1.0 & 0.3 &       \\
\hline 
\multicolumn{1}{l}{\textbf{RP3 \cite{DeledicquePapalexandris}:}} & 
\multicolumn{8}{c}{ $\gamma_s = 1.4, \quad \pi_s = 0, \quad \gamma_g = 1.4, \quad \pi_g = 0$}  \\ 
\hline
L & 1.0     & 0.9       & 2.5      & 1.0       & 0.0      &  1.0 & 0.9   & 0.10   \\
R & 1.0     & 0.0       & 1.0      & 1.2       & 1.0      &  2.0 & 0.2   &        \\
\hline 
\multicolumn{1}{l}{\textbf{RP4 \cite{Schwendeman}:}} & 
\multicolumn{8}{c}{ $\gamma_s = 3.0, \quad \pi_s = 3400, \quad \gamma_g = 1.35, \quad \pi_g = 0$}  \\
\hline
L & 1900.0   & 0.0   & 10.0    & 2.0 & 0.0   & 3.0 & 0.2 & 0.15   \\
R & 1950.0   & 0.0   & 1000.0  & 1.0 & 0.0   & 1.0 & 0.9 &       \\
\hline 
\multicolumn{1}{l}{\textbf{RP5 \cite{Schwendeman}:}} & 
\multicolumn{8}{c}{ $\gamma_s = 1.4, \quad \pi_s = 0, \quad \gamma_g = 1.4, \quad \pi_g = 0$}  \\
\hline 
L & 1.0    & 0.0   & 1.0  & 0.2 & 0.0   &  0.3 & 0.8 & 0.20  \\
R & 1.0    & 0.0   & 1.0  & 1.0 & 0.0   &  1.0 & 0.3 &      \\
\hline 
\multicolumn{1}{l}{\textbf{RP6 \cite{AndrianovWarnecke}:}} & 
\multicolumn{8}{c}{ $\gamma_s = 1.4, \quad \pi_s = 0, \quad \gamma_g = 1.4, \quad \pi_g = 0$}  \\
\hline 
L & 0.2068    & 1.4166   & 0.0416  & 0.5806 & 1.5833    &  1.375 & 0.1 & 0.10  \\
R & 2.2263    & 0.9366   & 6.0     & 0.4890 & -0.70138  &  0.986 & 0.2 &       \\
\hline
\end{tabular}
\end{center}
\label{tab.rpbn.ic}
\end{table}

\subsection{Explosion problems in multiple space dimensions} 
An explosion problem for the Baer-Nunziato equations can be solved both in two and in three space dimensions, 
after setting-up the following initial conditions
\begin{equation}
  \mathbf{u}(\mathbf x,0) = \left\{ \begin{array}{ccc} \mathbf{u}_i & \textnormal{ if } & r \leq R, \\ 
                                                       \mathbf{u}_o & \textnormal{ if } & r > R.
  \end{array} \right. 
\end{equation} 
where $\mathbf{x}$ is the vector of spatial coordinates, $r = \sqrt{\mathbf{x}^2}$,  
and $R$ is the radius of the initial discontinuity, which we have set equal to $0.5$ and to $0.4$ for the two-dimensional and for the three-dimensional explosion tests, respectively.
The inner and outer states $\mathbf{u}_i$ and $\mathbf{u}_o$ for the two test cases that we have considered are reported in detail in Table~\ref{tab.3d.explos.ic}. 
Analogously to the compressible Euler equations, a reference solution can be obtained after solving 
an equivalent  one-dimensional problem in cylindrical coordinates for the two-dimensional explosion, and 
in spherical coordinates for the three-dimensional one. 
This essentially implies that the equivalent one-dimensional problem contains additional algebraic source terms on the right hand side of the equations, which account for the use of curvilinear coordinates, 
see \cite{toro-book,USFORCE2}. The 1D reference solution has been computed using a classical second order TVD method with the Osher Riemann solver \cite{OsherUniversal,OsherNC}, using 5000 grid cells. 
For EPa the parameters for the refinement criterion are $\rho_{1,0}=\rho_{2,0}=1$ and for EPb we use $\rho_{1,0}=1000$ and $\rho_{2,0}=1$. 

\begin{table}[!b]
\caption{Inner and outer initial states for the two multidimensional explosion test problems. 
 } 
\renewcommand{\arraystretch}{1.0}
\begin{center}
\begin{tabular}{cccccccc|cccc}
\hline
 EPa   & $\rho_s$ & $p_s$ & $u_s$  & $\rho_g$ & $p_g$ & $u_g$ & $\phi_s$ & $\gamma_s$ & $\pi_s$ & $\gamma_g$ & $\pi_g$  \\
\hline
Inner  & 1.       & 1.  & 0.       & 0.5      & 1.0   &  0.   & 0.4      &  1.4       & 0.      & 1.4        & 0.       \\
Outer  & 2.       & 2.  & 0.       & 1.5      & 2.0   &  0.   & 0.8      &            &         &            &          \\
\hline
\hline
 EPb   & $\rho_s$ & $p_s$ & $u_s$  & $\rho_g$ & $p_g$ & $u_g$ & $\phi_s$ & $\gamma_s$ & $\pi_s$ & $\gamma_g$ & $\pi_g$  \\
\hline
Inner  & 800.    & 500.  & 0.   & 1.5      & 2.0   &  0.    & 0.4     &  3.0       & 100.    & 1.4        & 0.          \\
Outer  & 1000.   & 600.  & 0.   & 1.0      & 1.0   &  0.    & 0.3     &            &         &            &             \\
\hline
\end{tabular}
\end{center}
\label{tab.3d.explos.ic}
\end{table}
\subsubsection{2D computations}
\begin{figure}[!htbp] 
\begin{center}
\includegraphics[angle=0,width=6.6cm]{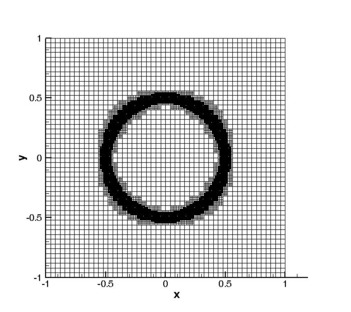}
\includegraphics[angle=0,width=6.6cm]{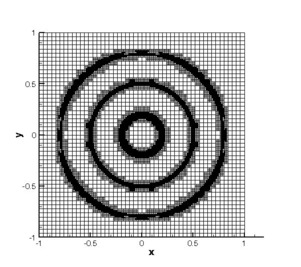}
\caption{  
AMR grid structure of the explosion test EPb in two space dimensions, at time $t=0$ (left panel)
and at the final time $t=0.2$ (right panel).
Two levels of refinement have been adopted ($\ell_{\max}=2$).
}
\label{fig_BNEP2Db-grid}
\end{center}
\end{figure}

\begin{figure}[!htbp] 
\begin{center}
\includegraphics[angle=0,width=6.6cm]{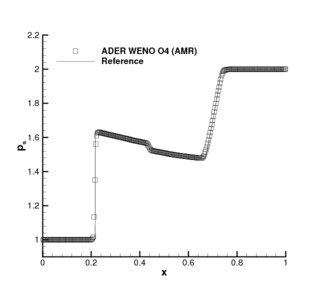}
\includegraphics[angle=0,width=6.6cm]{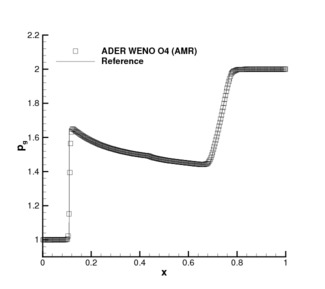}
\includegraphics[angle=0,width=6.6cm]{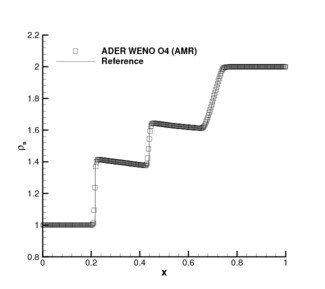}
\includegraphics[angle=0,width=6.6cm]{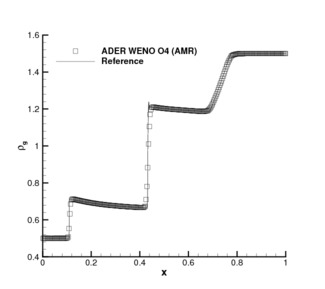}
\includegraphics[angle=0,width=6.6cm]{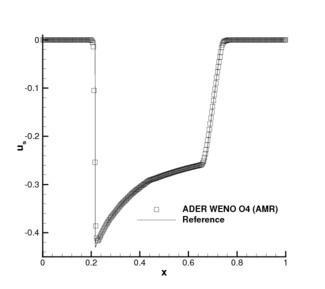}
\includegraphics[angle=0,width=6.6cm]{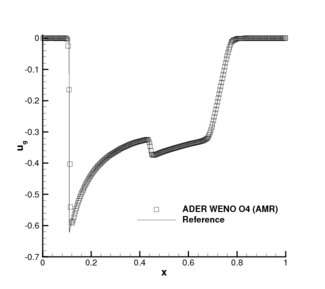}
\caption{  
Results of the 2D explosion problem EPa at time $t=0.2$. 
A cut of various quantities along the $x$-axis is reported, both for the solid (left panels) and for the gas phase (right panels).  The 1D reference solution is also shown for comparison.
Two levels of refinement have been adopted ($\ell_{\max}=2$).
}
\label{fig_BNEP2Da}
\end{center}
\end{figure}

\begin{figure}[!htbp] 
\begin{center}
\includegraphics[angle=0,width=6.6cm]{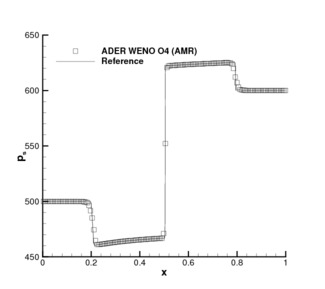}
\includegraphics[angle=0,width=6.6cm]{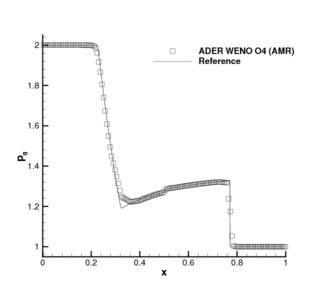}
\includegraphics[angle=0,width=6.6cm]{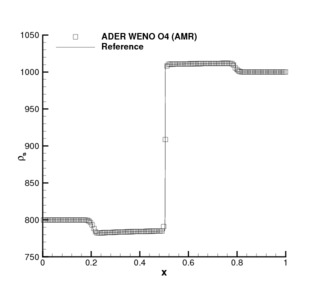}
\includegraphics[angle=0,width=6.6cm]{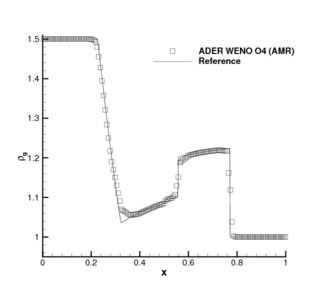}
\includegraphics[angle=0,width=6.6cm]{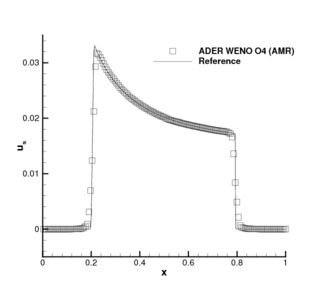}
\includegraphics[angle=0,width=6.6cm]{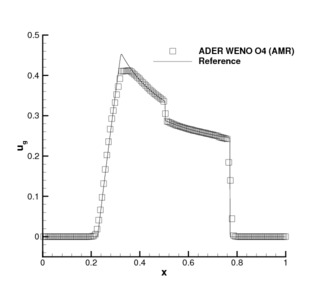}
\caption{  
Results of the 2D explosion problem EPb at time $t=0.2$. 
A cut of various quantities along the $x$-axis is reported, both for the solid (left panels) and for the gas phase (right panels).  The 1D reference solution is also shown for comparison.
Two levels of refinement have been adopted ($\ell_{\max}=2$).
}
\label{fig_BNEP2Db}
\end{center}
\end{figure}

We have first evolved the two models ${\rm EPa}$ and ${\rm EPb}$ in two spatial dimensions, by adopting a fourth order ADER-WENO-AMR scheme with two levels of refinement over a computational domain given by $[-1;1]\times[1;1]$. 
The level zero grid is composed by $50\times50$ cells, which are immediately refined according to the refinement criterion applied to the initial conditions. 
A representative example of the grid at time $t=0$ and at the final time $t=0.2$ is shown in Fig.~\ref{fig_BNEP2Db-grid} for the model ${\rm EPb}$. The final grid (right panel) is composed by $50020$ cells. Fig.~\ref{fig_BNEP2Da}  and Fig.~\ref{fig_BNEP2Db} report the results of the computation for the two models ${\rm EPa}$ and ${\rm EPb}$ by comparing them with the reference solution.
\subsubsection{3D computations}
\begin{figure}[!htbp] 
\begin{center}
\includegraphics[angle=0,width=6.6cm]{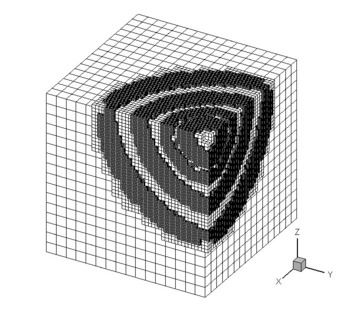}
\caption{  
AMR grid structure of the explosion test EPb in three space dimensions.
Two levels of refinement have been adopted ($\ell_{\max}=2$), starting from a uniform $34\times34\times34$ grid. 
}
\label{fig_BNEP3Db-grid}
\end{center}
\end{figure}

\begin{figure}[!htbp] 
\begin{center}
\includegraphics[angle=0,width=6.6cm]{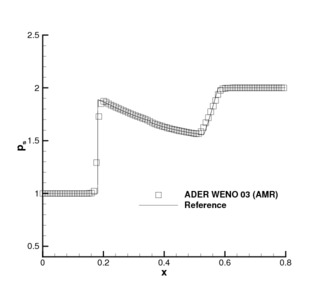}
\includegraphics[angle=0,width=6.6cm]{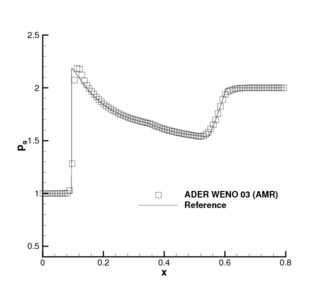}
\includegraphics[angle=0,width=6.6cm]{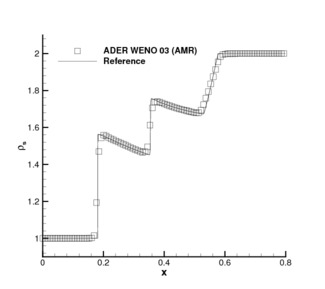}
\includegraphics[angle=0,width=6.6cm]{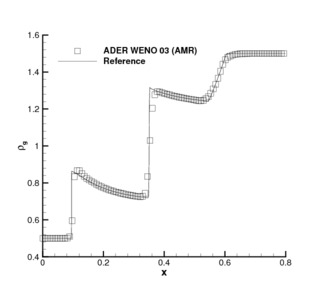}
\includegraphics[angle=0,width=6.6cm]{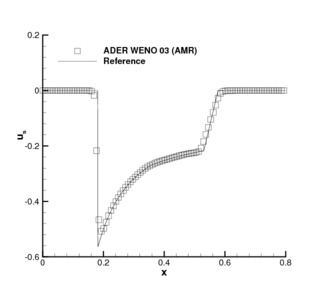}
\includegraphics[angle=0,width=6.6cm]{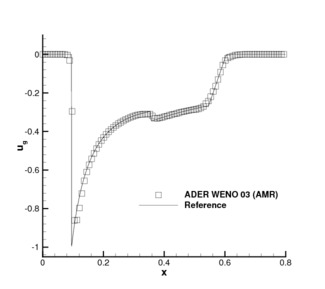}
\caption{  
Results of the 3D explosion problem EPa at time $t=0.15$. 
A cut of various quantities along the $x$-axis is reported, both for the solid (left panels) and for the gas phase (right panels).  The 1D reference solution is also shown for comparison.
Two levels of refinement have been adopted ($\ell_{\max}=2$), starting from a uniform $34\times34\times34$ grid. 
}
\label{fig_BNEP3Da}
\end{center}
\end{figure}

\begin{figure}[!htbp] 
\begin{center}
\includegraphics[angle=0,width=6.6cm]{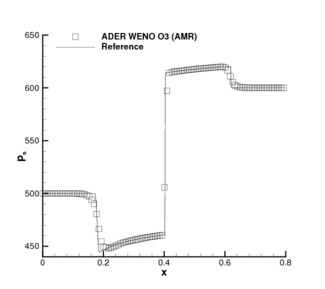}
\includegraphics[angle=0,width=6.6cm]{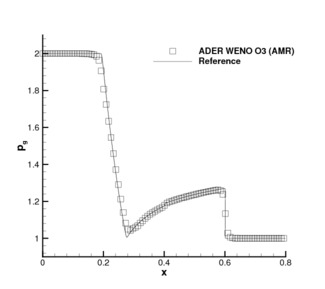}
\includegraphics[angle=0,width=6.6cm]{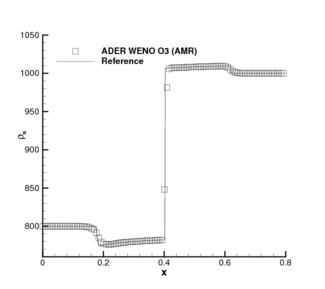}
\includegraphics[angle=0,width=6.6cm]{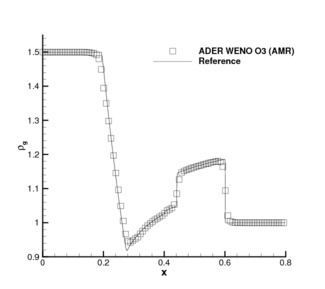}
\includegraphics[angle=0,width=6.6cm]{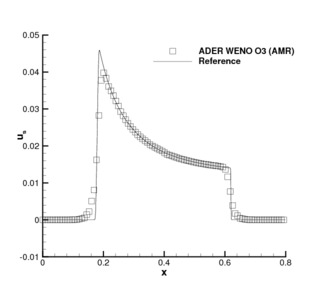}
\includegraphics[angle=0,width=6.6cm]{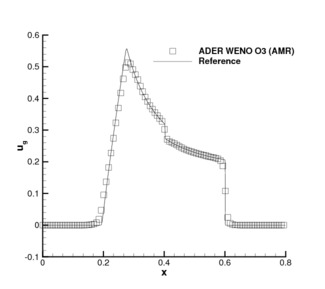}
\caption{  
Results of the 3D explosion problem EPb at time $t=0.15$. A cut of various quantities along the $x$-axis is reported, both for the solid (left panels) and for the gas phase (right panels).  The 1D reference solution is also shown for comparison.
Two levels of refinement have been adopted ($\ell_{\max}=2$), starting from a uniform $34\times34\times34$ grid. 
}
\label{fig_BNEP3Db}
\end{center}
\end{figure}

We have then evolved the same models ${\rm EPa}$ and ${\rm EPb}$ in three spatial dimensions, by adopting a third order ADER-WENO-AMR scheme with two levels of refinement over a computational domain given by  $[-1;1]\times[1;1]\times[-1;1]$. The level zero grid contains $34\times34\times34$ cells. The final grid at time $t=0.15$, shown as a representative example in Fig.~\ref{fig_BNEP3Db-grid}  
for the model ${\rm EPb}$, is composed by $3,833,016$ cells. Fig.~\ref{fig_BNEP3Da} and Fig.~\ref{fig_BNEP3Db}, on the other hand, report  the solution at time $t=0.15$, compared to the reference one. All the relevant features and waves of the solution are successfully resolved by the scheme, which remains essentially non-oscillatory and performs the correct grid refinement where this is needed.

\subsection{Shock-bubble interaction}

\begin{figure}[!htbp]
\begin{center}
\begin{tabular}{cc} 
\includegraphics[width=0.45\textwidth]{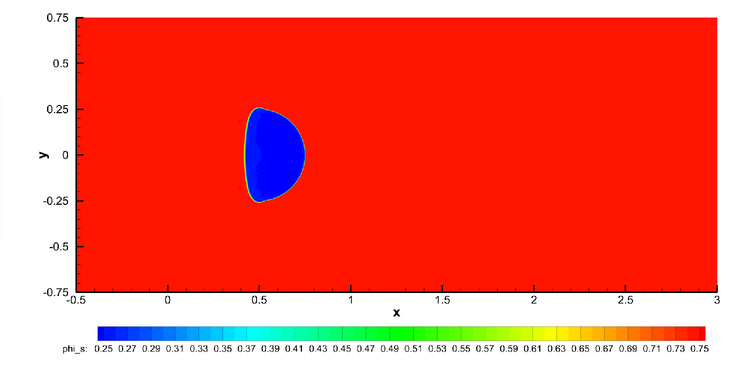}       & 
\includegraphics[width=0.45\textwidth]{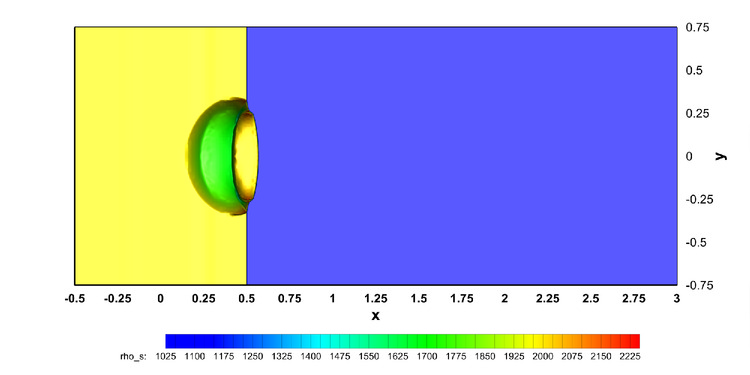}       \\ 
\includegraphics[width=0.45\textwidth]{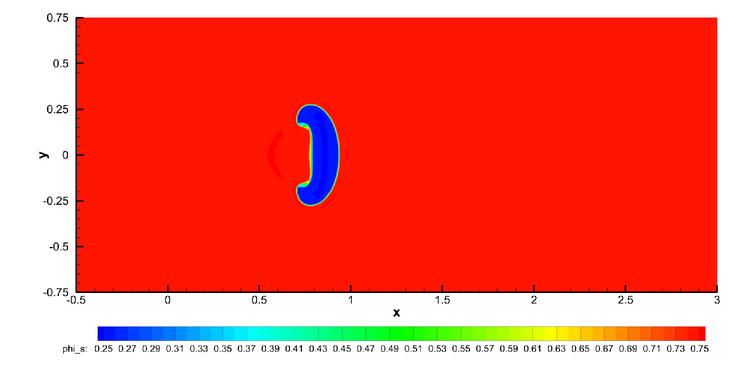}       & 
\includegraphics[width=0.45\textwidth]{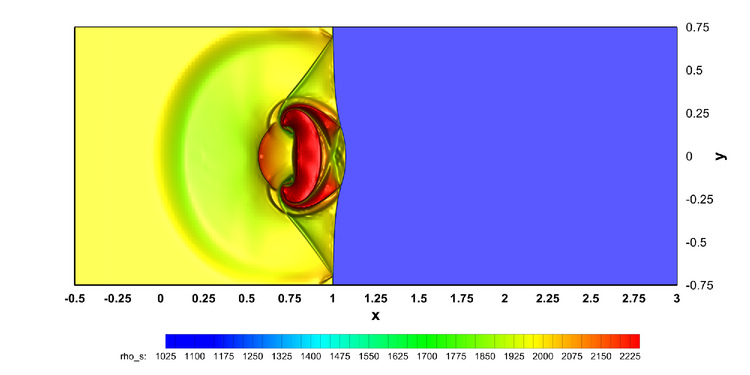}       \\ 
\includegraphics[width=0.45\textwidth]{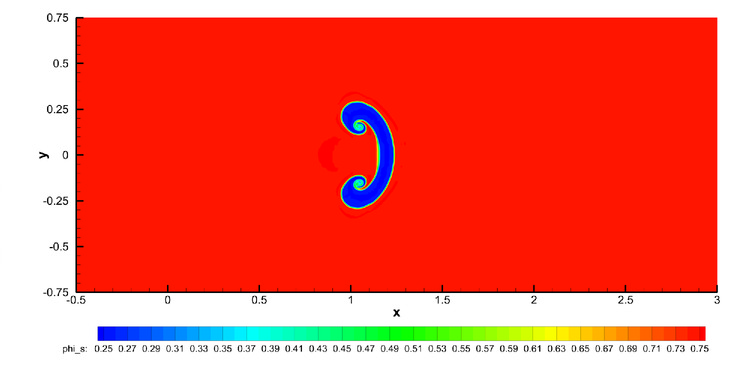}       & 
\includegraphics[width=0.45\textwidth]{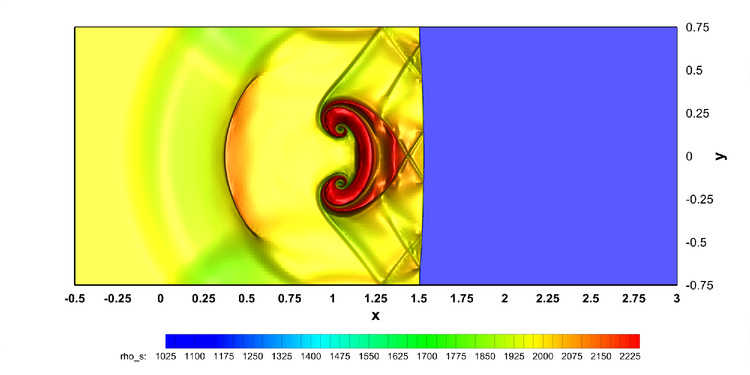}       \\ 
\includegraphics[width=0.45\textwidth]{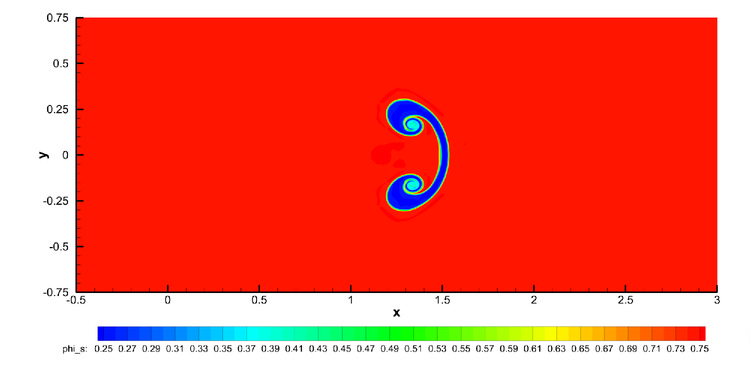}       & 
\includegraphics[width=0.45\textwidth]{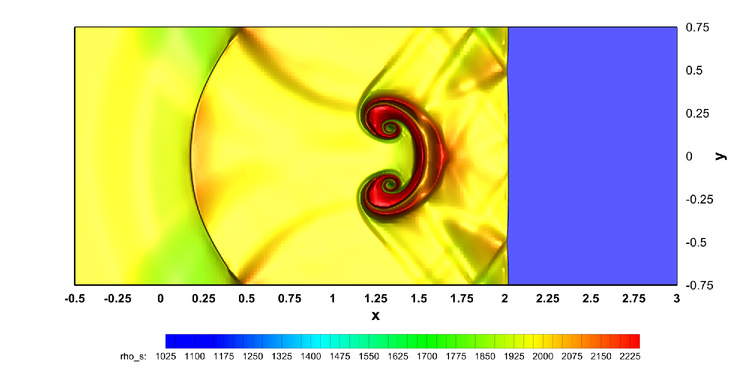}       \\ 
\includegraphics[width=0.45\textwidth]{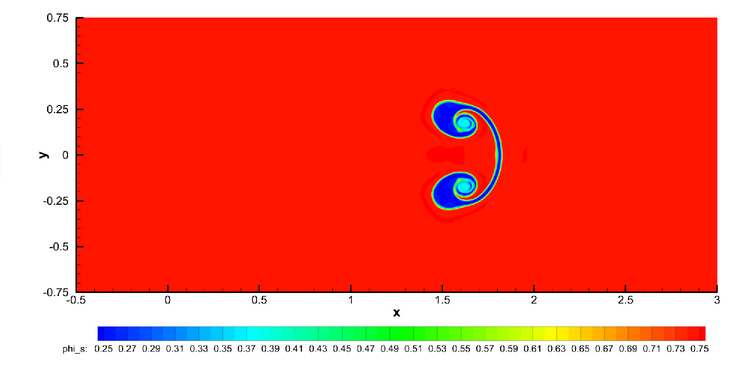}       & 
\includegraphics[width=0.45\textwidth]{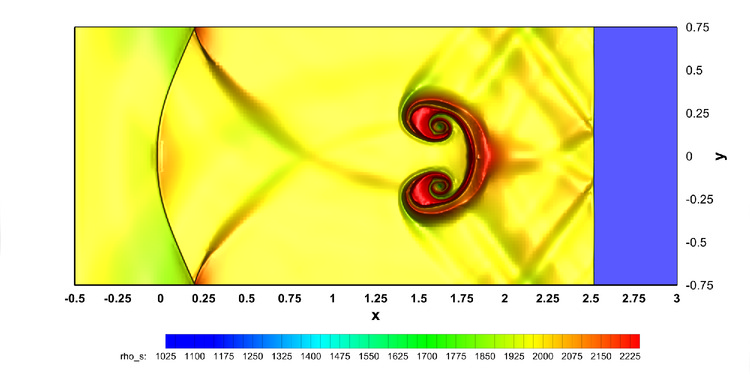}       \\ 
\end{tabular}
\caption{Liquid volume fraction (left) and liquid density (right) of the shock-bubble interaction problem at times $t=0.005$, $t=0.010$, $t=0.015$, $t=0.020$ and $t=0.025$ (from top to bottom) using a third order 
ADER-WENO scheme with AMR.} 
\label{fig.bubble}
\end{center}
\end{figure}

\begin{figure}[!htbp]
\begin{center}
\includegraphics[width=0.65\textwidth]{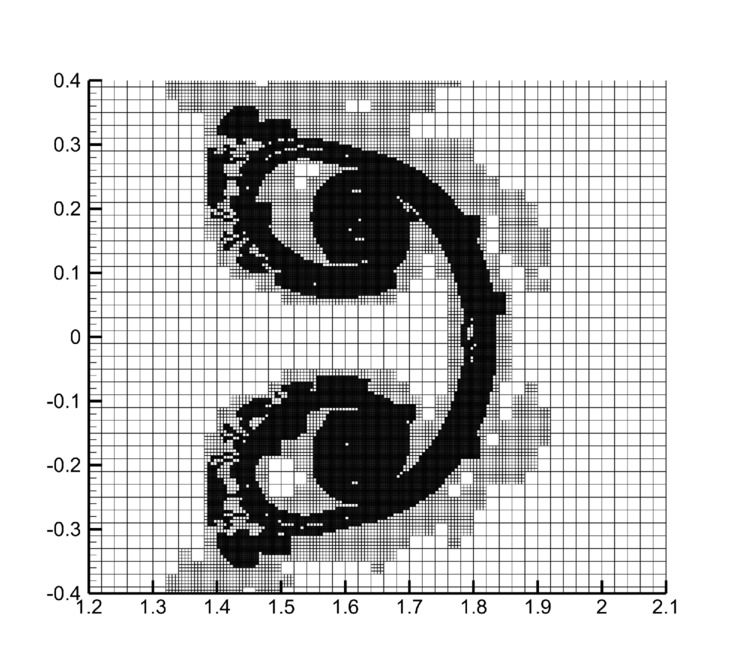}        
\caption{Zoom into the AMR grid at time $t=0.025$ for the shock-bubble interaction problem.}
\label{fig.bubble.grid}
\end{center}
\end{figure}

In this section a strongly simplified shock-bubble interaction-type problem is considered. The initial condition is given by a Riemann problem that leads to a strong isolated shock wave travelling 
with a shock speed of $s=100$. The pressure jumps over six orders of magnitude across the shock. The initial discontinuity is located at $x=0$ and the left and right initial 
states, which are connected by a Hugoniot curve, are summarized in Table \ref{tab.bnbubble}. A bubble of radius $R=0.25$ is located at $x=0.5$, $y=0$. The state inside the 
bubble is also given in Table \ref{tab.bnbubble}. 
The parameters for the equation of state of each phase are $\gamma_1 = 3.0$, $\pi_1 = 100$, $\gamma_2 = 1.4$ and $\pi_2 = 0$. The problem is solved until time 
$t=0.0025$ in a computational domain $\Omega = [-0.5;3.0] \times [-0.75;0.75]$ that is discretized on the level $\ell=0$ grid using $70 \times 30$ grid zones. 
Two levels of refinement are used ($\ell_{\max}=2$) with a refinement factor of $\mathfrak{r}=4$, which corresponds to a uniform fine grid resolution of $ 1120 \times 480$.  
Periodic boundary conditions are employed in $y$ direction and Dirichlet boundary conditions corresponding to the left and right state are imposed in $x$ direction. 
We choose $\rho_{1,0}=1000$ and $\rho_{2,0}=1$ for the indicator function $\Phi$. 
The evolution of the liquid density $\rho_1$ and the liquid volume fraction $\phi_1$ are depicted for various times in Fig. \ref{fig.bubble}. One can observe how the 
bubble is compressed and accelerated by the incident shock wave and how Richtmyer-Meshkov-type flow instabilities occur at the bubble border at later times. One 
further notes the shock wave reflected by the bubble. A zoom into the AMR grid at the final time $t=0.0025$ is depicted in Fig. \ref{fig.bubble.grid}. The results of this 
section are only considered as \textit{qualitative}, in order to test the robustness of our high order AMR method in the presence of strong shock waves in liquids. 
Here, the high order ADER-WENO AMR scheme has been applied to the full seven equation Baer-Nunziato model, but without taking into account any interphase drag and pressure relaxation, hence
no comparison with experimental data can be made for our results. A very detailed \textit{quantitative} study of shock-bubble interactions with comparison against experimental data has been
carried out in \cite{QuirkKarni} using second order accurate high resolution shock capturing schemes together with AMR. 
\begin{table}[!b]
\caption{Left (L), right (R) and bubble (B) state for the shock-bubble interaction problem. } 
\renewcommand{\arraystretch}{1.0}
\begin{center}
\begin{tabular}{lccccccc}
\hline
       & $\rho_1$ & $u_1$ & $p_1$  & $\rho_2$ & $u_2$ & $p_2$ & $\phi_1$    \\
\hline
Left (L)    & 1999.939402       & 49.998485    & 4999849.5      & 1.0      & 0.0   &  1.0   & 0.75        \\
Right (R)   & 1000.0            & 0.0          & 1.0            & 1.0      & 0.0   &  1.0   & 0.75        \\
Bubble (B)  & 1000.0            & 0.0          & 1.0            & 1.0      & 0.0   &  1.0   & 0.25        \\ 
\hline
\end{tabular}
\end{center}
\label{tab.bnbubble}
\end{table}

\subsection{An application to complex free surface flows} 
\label{sec.db2d}

A reduced version of the Baer--Nunziato model (\ref{ec.BN}) can also be used for the simulation of complex non--hydrostatic free surface flows, as proposed in \cite{DIM2D,DIM3D}. The reduced three--equation version of the PDE  system (\ref{ec.BN}) has been obtained in \cite{DIM2D,DIM3D} after introducing the following simplifying assumptions:  i) all pressures are relative with respect to the atmospheric reference pressure, which is assumed to be  constant and zero everywhere and for all times; 
ii) the influence of the gas phase onto the liquid can be neglected, hence the evolution equations for the gas phase can be dropped in (\ref{ec.BN});    
iii) the pressure of the liquid phase is computed by the Tait equation of state \cite{Batchelor}, which is a barotropic EOS and therefore also allows the energy equation of the liquid to be removed from the governing PDE system. 
The Tait EOS reads 
	\begin{equation}
  p_1 = k_0 \left[\left(\frac{\rho_1}{\rho_0}\right)^\gamma-1\right],
  \label{eq:TaitEOS}
  \end{equation}
where $\rho_1,\rho_0$ are the liquid density and the reference liquid density at standard conditions, respectively. The compressibility of the fluid is governed by the constants $k_0$ and $\gamma$. 
To avoid low Mach number problems, an artificial Mach number of the order $M \in [0.1; 0.3]$ is chosen for typical environmental free surface problems. If real compressibility effects play a role,
such as in industrial high pressure applications, the proper values $k_0=3.2 \cdot 10^8$ $Pa$, $\rho_0=1000$ $kg/m^3$ and $\gamma=7$ for real water can be chosen, which lead to the real sound speed 
in water of approximately $c = 1500$ m/s. 

Introducing the above-mentioned simplifications into system (\ref{ec.BN}) yields the following reduced three-equation model \cite{DIM2D,DIM3D}:  
\begin{equation}
\renewcommand{\arraystretch}{1.5}
\begin{array}{ll}
&\frac{\partial}{\partial t}\left(\phi \rho \right)+\nabla \cdot \left(\phi \rho \mathbf{v}\right) = 0, \\
&\frac{\partial}{\partial t}\left(\phi \rho \mathbf{v}\right)+\nabla \cdot \left(\phi\left(\rho\mathbf{v}\mathbf{v} + p \mathbf{I} \right)\right) = \phi \rho \mathbf{g}, \\
&\frac{\partial}{\partial t}\phi + \mathbf{v}\cdot\nabla\phi = 0.  
\label{eq:ViscreducedBN}
\end{array}
\end{equation}
The subscript $1$ has been dropped to ease notation. The mass and momentum equations are fully conservative in the system above since the interface pressure $p_I=p_2=0$, while the advection equation for 
the volume fraction remains non--conservative. In (\ref{eq:ViscreducedBN}) the state vector is $\mathbf{Q}=\left(\phi\rho,\phi\rho\mathbf{v},\phi\right)$ and $\mathbf{g}$ is the gravity vector 
acting along the vertical direction $y$, i.e. $\mathbf{g} = (0,-g)$ with $g=9.81$ $m/s^2$. Further details and thorough validations of this model can be found in  \cite{DIM2D,DIM3D}. 

\begin{figure}[!htbp]
\begin{center}
\begin{tabular}{cc} 
\includegraphics[width=0.45\textwidth]{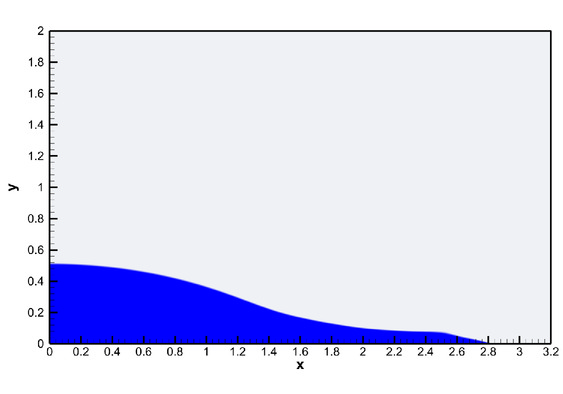}       & 
\includegraphics[width=0.45\textwidth]{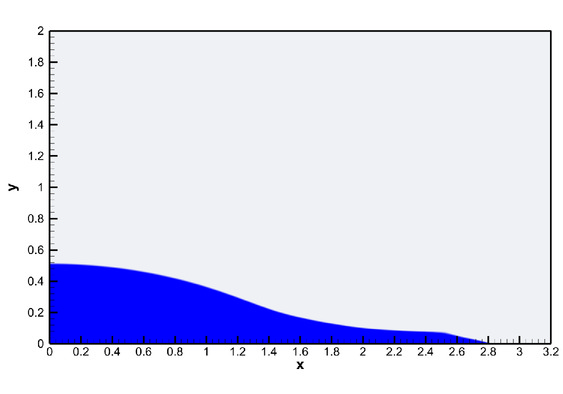}       \\ 
\includegraphics[width=0.45\textwidth]{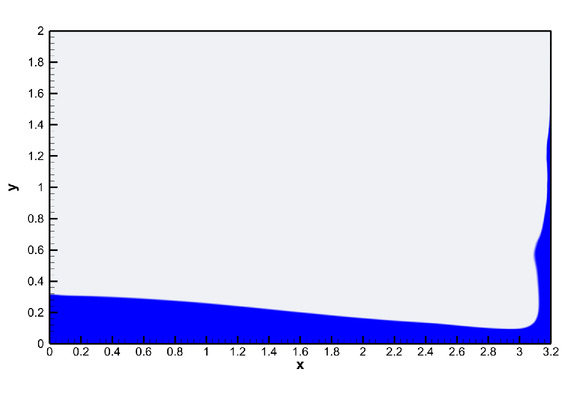}       & 
\includegraphics[width=0.45\textwidth]{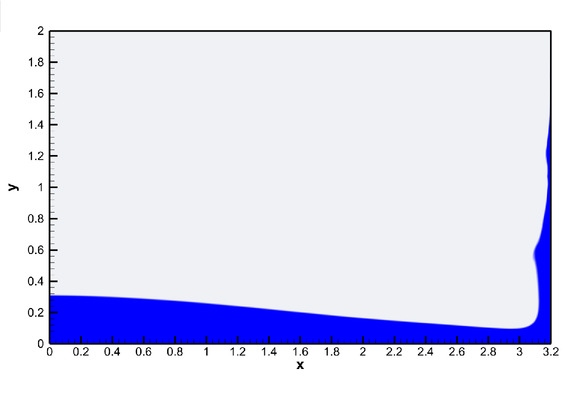}       \\ 
\includegraphics[width=0.45\textwidth]{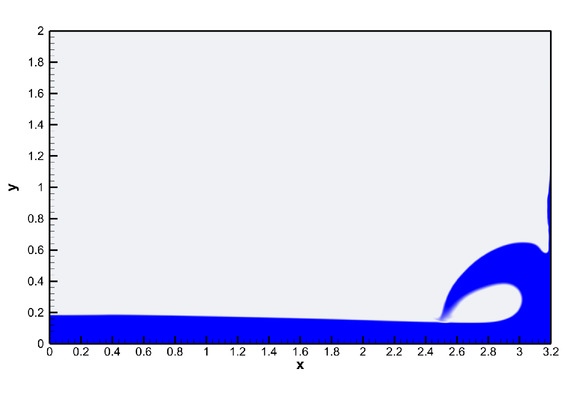}       & 
\includegraphics[width=0.45\textwidth]{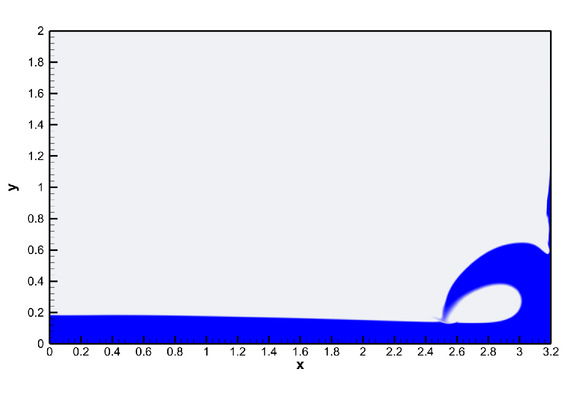}       
\end{tabular}
\caption{Dambreak and wave-impact problem at times $t=0.5$, $t=1.0$ and $t=1.5$ using a third order ADER-WENO scheme with AMR (left) and with a uniform fine grid (right).}
\label{fig.db2d}
\end{center}
\end{figure}

\begin{figure}[!htbp]
\begin{center}
\includegraphics[width=0.75\textwidth]{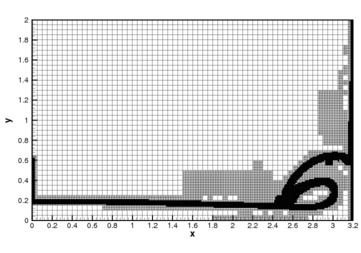}        
\caption{AMR grid for the dambreak and wave impact problem at time $t=1.5$.}
\label{fig.db2d.grid}
\end{center}
\end{figure}

\begin{figure}[!]
\begin{center}
\includegraphics[width=0.65\textwidth]{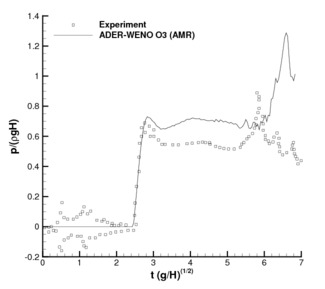}        
\caption{Pressure evolution at the wall. Comparison of experimental data with 
the computational results.}
\label{fig.db2d.p}
\end{center}
\end{figure}

The last test problem of this paper where we apply the high order ADER-WENO scheme with adaptive mesh refinement consists of a dambreak with successive wave impact against a vertical wall. 
In this test problem the free surface is strongly deformed and the impact against the wall also leads to wave breaking. The setup of this test case is taken from \cite{Colagrossi,SPH3D}.  
The rectangular computational domain is defined as $\Omega = [0;3.2] \times [0;2.0]$ with reflective wall boundaries on all sides of the domain apart from the top, which is assumed to be transmissive. 
The initial liquid domain is $\Omega_l = [0;1.2] \times [0;0.6]$, where $\phi = 0.999$ is set. The initial velocity is zero and the initial density distribution is chosen such as to obtain a hydrostatic
pressure distribution inside the liquid, see \cite{DIM2D} for details. Outside the liquid domain, we set the pressure to zero and a low but finite value for the liquid volume fraction is chosen 
($\phi = 10^{-3}$). We use $g=9.81$ and $k_0 = 2.0 \cdot 10^5$. The level zero grid contains only $60 \times 40$ cells. We use two levels of refinement $\ell_{\max}=2$ and a refinement factor of 
$\mathfrak{r}=4$. This corresponds to an equivalent resolution on a uniform fine mesh of $960 \times 640$ elements. A third order ADER-WENO method is used together with the little dissipative Osher scheme. 
The contours of the liquid volume fraction $\phi$ are depicted for three output times in Fig. \ref{fig.db2d}. The results of the present high order AMR method are compared with the results on a uniform fine 
grid. According to Fig. \ref{fig.db2d} the two solutions match very well. The AMR grid at the final time $t=1.5$ is depicted in Fig. \ref{fig.db2d.grid} and contains 47710 cells, while the uniform fine grid 
contains 614400 cells. This means that for the present test problem the use of an AMR technique allows to obtain almost the same results as with a uniform fine grid, but with 12.88 times less elements. 
Concerning the AMR overhead we also report the average CPU time per element update (EU) in Table \ref{tab.db2d.compare}, which is the total wallclock time divided by the total number of element updates, normalized
with respect to the uniform grid. This figure indicates the overhead introduced by the AMR machinery. For the present test problem, using a third order ADER-WENO AMR scheme it is only 19 \%, which agrees with
the data for the AMR overhead published in \cite{Dumbser2012a}.  
One can note how the AMR algorithm refines the grid in the vicinity of the free surface, which is resolved in a very sharp manner. The flow evolution computed with the present method is also in good agreement 
with the 3D SPH computations performed in \cite{SPH3D} and with the results presented in \cite{DIM2D,DIM3D}. In Fig. \ref{fig.db2d.p} we compare the experimental data for the pressure at the wall \cite{ExpCol} 
with the numerical results obtained with the method proposed in this paper. The amplitude and the arrival time of the first pressure peak due to the wave impact at the wall is captured correctly. 

\begin{table}[!b]   
\caption{Memory and CPU time comparison of the third order ADER-WENO AMR method and ADER-WENO on a uniform fine grid for the dambreak and wave impact problem of Section \ref{sec.db2d}. 
         Memory consumption is measured in number of real elements and CPU time is normalized with respect to the wallclock time for the fine uniform mesh. The normalized average time per element update (EU) 
         is given in the last row to quantify the overhead of the AMR framework.} 
\begin{center} 
\renewcommand{\arraystretch}{1.0}
\begin{tabular}{cccc} 
\hline
                           & AMR     & Uniform &  ratio \\ 
\hline
Number of real cells       & 47710   & 614400  & 12.877 \\ 
Total CPU time						 & 0.1134  & 1.0     & 8.816  \\ 
Average time per EU        & 1.1923  & 1.0     & 0.839  \\  
\hline
\end{tabular} 
\end{center}
\label{tab.db2d.compare}
\end{table}

\section{Conclusions}
\label{sec:concl}

In the present paper the first better than second order accurate path-conservative one-step WENO finite volume 
scheme with adaptive mesh refinement (AMR) has been presented for the solution of non-conservative hyperbolic
systems. The method has been applied in particular to the Baer-Nunziato model of compressible multiphase flows. 
It has been shown via a numerical convergence study that the proposed numerical approach reaches the designed high 
order of accuracy in space and time. Furthermore, the accuracy and robustness of the method has been validated on a 
large set of test problems, for which either exact or other reference solutions were available. It has been clearly 
shown that the use of AMR can lead to a significant reduction in grid elements and CPU time compared to simulations 
performed on fine uniform grids. Compared to our previous publications \cite{USFORCE2,OsherNC}, where high order 
schemes have been employed without AMR, the material interfaces are much better resolved using the present high order 
AMR technique. This becomes particularly evident in the results obtained for the shock tube problems solved in Section 
\ref{sec.bnrp}. 

Future applications of the present high order one-step AMR methodology may concern the simulation of bubbles with phase
transition \cite{DreyerHantke} and chemically reacting compressible multiphase flows. For the latter it has already been 
shown before in literature that AMR techniques can be very useful to resolve all the length scales involved in such 
multi-scale problems \cite{HelzelLevequeWarnecke}. 
Other relevant mathematical models that could be treated in the framework outlined in this paper include the nonconservative 
debris flow model of Pitman and Le \cite{PitmanLe} as well as single and multi-layer shallow water equations 
\cite{Castro2006,Pares2006,CastroCendon,AbgrallKarni3}. In particular for multi-scale tsunami wave simulations the combination of a high 
order method that allows a coarse grid to discretize the wave propagation in the ocean with a fine grid that allows to follow 
and resolve all the relevant details on the coastline may be useful. 

\section*{Acknowledgments}

The research conducted here has been financed by the European Research Council under 
the European Union's Seventh Framework Programme (FP7/2007-2013) in the frame of the research 
project \textit{STiMulUs}, ERC Grant agreement no. 278267. A.H. thanks Fundaci\'on Caja Madrid 
(Spain) for its financial support by a grant under the programme \textit{Becas de movilidad 
para profesores de las universidades p\'ublicas de Madrid}.






\bibliography{AMR-Nonconservative}
\bibliographystyle{plain}







\end{document}